%%%% prevent double loading:
\expandafter\ifx\csname mthreemacsloaded\endcsname\relax\else \fi

\magnification1100
\input amstex

%%% Hack of Plain TeX correction and style macros 
%%% written by Walter Neumann and Larry Siebenmann:

 \catcode`\@=11
 \let\wlog@ld\wlog
 \def\wlog#1{\relax}

 \newif\ifIN@
 \def\m@rker{\m@@rker}
 \def\IN@{\expandafter\INN@\expandafter}
 \long\def\INN@0#1@#2@{\long\def\NI@##1#1##2##3\ENDNI@
    {\ifx\m@rker##2\IN@false\else\IN@true\fi}%
     \expandafter\NI@#2@@#1\m@rker\ENDNI@}
  \newtoks\Initialtoks@  \newtoks\Terminaltoks@
  \def\SPLIT@{\expandafter\SPLITT@\expandafter}
  \def\SPLITT@0#1@#2@{\def\TTILPS@##1#1##2@{%
     \Initialtoks@{##1}\Terminaltoks@{##2}}\expandafter\TTILPS@#2@}
  \newtoks\Trimtoks@

 \def\ForeTrim@{\expandafter\ForeTrim@@\expandafter}
 \def\ForePrim@0 #1@{\Trimtoks@{#1}}
 \def\ForeTrim@@0#1@{\IN@0\m@rker. @\m@rker.#1@%
     \ifIN@\ForePrim@0#1@%
     \else\Trimtoks@\expandafter{#1}\fi}
 
  \def\Trim@0#1@{%
      \ForeTrim@0#1@%
      \IN@0 @\the\Trimtoks@ @%
        \ifIN@
             \SPLIT@0 @\the\Trimtoks@ @\Trimtoks@\Initialtoks@
             \IN@0\the\Terminaltoks@ @ @%
                 \ifIN@
                 \else \Trimtoks@ {FigNameWithSpace}%
                 \fi
        \fi
      }

  %%% Math Bolds
  \font\titlebold=cmbx12 scaled 1200
  \font\twelvebold=cmbx12
  \font\tenbold=cmbx10
  \font\ninebold=cmbx9
  \font\sevenbold=cmbx7
  \font\fivebold=cmbx5

  \input amssym.def \input amssym
  %%% point sizes not loaded by amssym.def:
     \font\titlemsa=msam10 at 14.4pt
     \font\titlemsb=msbm10 at 14.4pt
     \font\titleeufm=eufm10 at 14.4pt
     \font\twelvemsa=msam10 scaled 1200
     \font\twelvemsb=msbm10 scaled 1200
     \font\twelveeufm=eufm10 scaled 1200
     \font\ninemsa=msam9
     \font\ninemsb=msbm9
     \font\nineeufm=eufm9

   %%% Cyrillic fonts (for accents and input, see ams cyr doc)
   \ifx\cyrfam\undefined
   \else
     \immediate\write16{}%
     \message{ !!! cyr fonts already defined. !!! }
     \message{ --- edit out superfluous font defs? }
   \fi
   \newfam\cyrfam
       \font\titlecyr=wncyr10 scaled 1440 %%% no caps?
       \font\twelvecyr=wncyr10 scaled 1200
       \font\tencyr=wncyr10
       \font\ninecyr=wncyr9
       \font\sevencyr=wncyr7
       \font\sixcyr=wncyr6

   %%% Euler script fonts (replacing caligraphic):
   \newfam\eusmfam
       \font\titleeusm=eusm10 scaled 1440
       \font\twelveeusm=eusm10 scaled 1200
       \font\teneusm=eusm10
       \font\nineeusm=eusm9
       \font\seveneusm=eusm7
       
       \font\fiveeusm=eusm5  

    \def\cal#1{{\fam\eusmfam\relax#1}}
    
\let\Cal\cal

 %%% Some fonts not loaded by plain
    \font\ninemrm=cmr9 %% new name for 9 pt math roman
    \font\ninei=cmmi9
    \font\ninesy=cmsy9 
    \skewchar\ninei='177
    \skewchar\ninesy='60

  \font\twelvemrm=cmr10 at 12pt %% new name
  \font\twelvei=cmmi10 at 12pt
  \font\twelvesy=cmsy10 at 12pt
 % \font\twelveex=cmex10 at 12pt

  \font\titlemrm=cmr10 at 14.4pt %% new name
  \font\titlei=cmmi10 at 14.4pt
  \font\titlesy=cmsy10 at 14.4pt
 % \font\titleex=cmex10 at 14.4pt

 %%%% Miscellanious font definitions

  \def\Smallfonts{\ninepoint}

  \def\Hfont{\titlepoint\bf}
  \def\Authorfont{\twelvepoint\it}
  \def\HHfont{\twelvepoint\bf}
  \def\HHHfont{\bf}
  % automatically smaller in 9 point parts
  \def\Bibfont{\tenbf}
  \def\Coordfont{\nineit }% defined in osuPSfnt.sty

  \def \thfont {\bf }
  \def \pffont {\it\itSpacing }
  \def \rkfont {\bf }
  \def \dffont {\bf }
  \def \egfont {\bf }

 %%%%% NINEPOINT %%%%%
 \def\ninepoint{%
  \def\rm{\fam0\ninerm}%
    \textfont0=\ninemrm  \scriptfont0=\sevenrm  \scriptscriptfont0=\fiverm
    \textfont1=\ninei    \scriptfont1=\seveni   \scriptscriptfont1=\fivei
  \def\mit{\fam1\ninei}%
  \def\oldstyle{\fam1\ninei}%
    \textfont2=\ninesy   \scriptfont2=\sevensy  \scriptscriptfont2=\fivesy
    \textfont3=\tenex    \scriptfont3=\tenex    \scriptscriptfont3=\tenex
  \def\it{\fam\itfam\nineit}%
    \textfont\itfam=\nineit
  \def\bf{\ifmmode\fam\bffam\else\ninebf\fi}%
    \textfont\bffam=\ninebold 
    \scriptfont\bffam=\sevenbold 
    \scriptscriptfont\bffam=\fivebold%
  \def\msa{\fam\msafam\ninemsa}%
    \textfont\msafam=\ninemsa 
    \scriptfont\msafam=\sevenmsa
    \scriptscriptfont\msafam=\fivemsa%
  \def\msb{\fam\msbfam\ninemsb}%
    \textfont\msbfam=\ninemsb%
    \scriptfont\msbfam=\sevenmsb%
    \scriptscriptfont\msbfam=\fivemsb%
  \def\eufm{\fam\eufmfam\nineeufm}%
    \textfont\eufmfam=\nineeufm
    \scriptfont\eufmfam=\seveneufm
    \scriptscriptfont\eufmfam=\fiveeufm
   \def\eusm{\fam\eusmfam\nineeusm}%
     \textfont\eusmfam=\nineeusm
     \scriptfont\eusmfam=\seveneusm
     \scriptscriptfont\eusmfam=\fiveeusm
   \def\cyr{\fam\cyrfam\ninecyr}%
     \textfont\cyrfam=\ninecyr
     \scriptfont\cyrfam=\sevencyr
     \scriptscriptfont\cyrfam=\sixcyr%%
  \setbox\strutbox=\hbox{\vrule
      height7pt depth3pt width0pt}%
   \baselineskip=10.8pt\rm}

 \let\eightpoint\ninepoint % we do not use eightpoint

 %%%%% FONTS AT TENPOINT %%%%%
 \def\tenpoint{%
  \def\rm{\fam0\tenrm}%
    \textfont0=\tenmrm \scriptfont0=\sevenrm \scriptscriptfont0=\fiverm%
  \def\mit{\fam1\teni}%
  \def\oldstyle{\fam1\teni}%
    \textfont1=\teni   \scriptfont1=\seveni  \scriptscriptfont1=\fivei%
    \textfont2=\tensy  \scriptfont2=\sevensy \scriptscriptfont2=\fivesy%
    \textfont3=\tenex  \scriptfont3=\tenex   \scriptscriptfont3=\tenex%
  \def\it{\fam\itfam\tenit}%
    \textfont\itfam=\tenit%
  \def\bf{\ifmmode\fam\bffam\else\tenbf\fi}%
    \textfont\bffam=\tenbold% was tenbold for osu
    \scriptfont\bffam=\sevenbold%
    \scriptscriptfont\bffam=\fivebold%
  \def\msa{\fam\msafam\tenmsa}%
    \textfont\msafam=\tenmsa%
    \scriptfont\msafam=\sevenmsa%
    \scriptscriptfont\msafam=\fivemsa%
  \def\msb{\fam\msbfam\tenmsb}%
    \textfont\msbfam=\tenmsb%
    \scriptfont\msbfam=\sevenmsb%
    \scriptscriptfont\msbfam=\fivemsb%
  \def\eufm{\fam\eufmfam\teneufm}%
   \textfont\eufmfam=\teneufm
   \scriptfont\eufmfam=\seveneufm
   \scriptscriptfont\eufmfam=\fiveeufm
   \def\eusm{\fam\eusmfam\teneusm}%
    \textfont\eusmfam=\teneusm
    \scriptfont\eusmfam=\seveneusm
    \scriptscriptfont\eusmfam=\fiveeusm
   \def\cyr{\fam\cyrfam\tencyr}%
    \textfont\cyrfam=\tencyr
    \scriptfont\cyrfam=\sevencyr
    \scriptscriptfont\cyrfam=\sixcyr%%
  \setbox\strutbox=\hbox{\vrule %
      height8.5pt depth3.5ptwidth0pt}%
  \baselineskip=\StdBaselineskip\rm}

 %%%%% FONTS AT TWELVEPOINT %%%%%
 \def\twelvepoint{%
  \def\rm{\fam0\twelverm}%
    \textfont0=\twelvemrm \scriptfont0=\tenmrm \scriptscriptfont0=\sevenrm
    \textfont1=\twelvei   \scriptfont1=\teni   \scriptscriptfont1=\seveni
  \def\mit{\fam1\twelvei}%
  \def\oldstyle{\fam1\twelvei}%
    \textfont2=\twelvesy  \scriptfont2=\tensy  \scriptscriptfont2=\sevensy
    \textfont3=\tenex  \scriptfont3=\tenex  \scriptscriptfont3=\tenex
  \def\it{\fam\itfam\twelveit}%
    \textfont\itfam=\twelveit
  \def\bf{\ifmmode\fam\bffam\else\twelvebf\fi}%
    \textfont\bffam=\twelvebold
    \scriptfont\bffam=\tenbold%
    \scriptscriptfont\bffam=\sevenbold%
  \def\msa{\fam\msafam\twelvemsa}%
    \textfont\msafam=\twelvemsa%
    \scriptfont\msafam=\tenmsa%
    \scriptscriptfont\msafam=\sevenmsa%
  \def\msb{\fam\msbfam\twelvemsb}%
    \textfont\msbfam=\twelvemsb%
    \scriptfont\msbfam=\tenmsb%
    \scriptscriptfont\msbfam=\sevenmsb%
  \def\eufm{\fam\eufmfam\twelveeufm}%
   \textfont\eufmfam=\twelveeufm
   \scriptfont\eufmfam=\teneufm
   \scriptscriptfont\eufmfam=\seveneufm
   \def\eusm{\fam\eusmfam\twelveeusm}%
    \textfont\eusmfam=\twelveeusm
    \scriptfont\eusmfam=\teneusm
    \scriptscriptfont\eusmfam=\seveneusm
   \def\cyr{\fam\cyrfam\tencyr}%
    \textfont\cyrfam=\twelvecyr
    \scriptfont\cyrfam=\tencyr
    \scriptscriptfont\cyrfam=\sevencyr%%
  \setbox\strutbox=\hbox{\vrule
      height10.2pt depth4.55pt width0pt}%
  \baselineskip=14pt\rm}

 %%%%% FONTS AT TITLEPOINT %%%%%
 \def\titlepoint{%
    \textfont0=\titlemrm \scriptfont0=\twelvemrm \scriptscriptfont0=\tenmrm
    \textfont1=\titlei   \scriptfont1=\twelvei   \scriptscriptfont1=\teni
  \def\mit{\fam1\titlei}%
  \def\oldstyle{\fam1\titlei}%
    \textfont2=\titlesy  \scriptfont2=\twelvesy  \scriptscriptfont2=\tensy
    \textfont3=\tenex% math ext not avail in varying sizes??
    \scriptfont3=\tenex
    \scriptscriptfont3=\tenex
  \def\it{\fam\itfam\titleit}%
    \textfont\itfam=\titleit
  \def\bf{\ifmmode\fam\bffam\else\titlebf\fi}%
    \textfont\bffam=\titlebold
    \scriptfont\bffam=\twelvebold%
    \scriptscriptfont\bffam=\tenbold%
  \def\msa{\fam\msafam\titlemsa}%
    \textfont\msafam=\titlemsa%
    \scriptfont\msafam=\twelvemsa%
    \scriptscriptfont\msafam=\tenmsa%
  \def\msb{\fam\msbfam\titlemsb}%
    \textfont\msbfam=\titlemsb%
    \scriptfont\msbfam=\twelvemsb%
    \scriptscriptfont\msbfam=\tenmsb%
  \def\eufm{\fam\eufmfam\titleeufm}%
    \textfont\eufmfam=\titleeufm
    \scriptfont\eufmfam=\twelveeufm
    \scriptscriptfont\eufmfam=\teneufm
   \def\eusm{\fam\eusmfam\titleeusm}%
     \textfont\eusmfam=\titleeusm
     \scriptfont\eusmfam=\twelveeusm
     \scriptscriptfont\eusmfam=\teneusm
   \def\cyr{\fam\cyrfam\tencyr}%
    \textfont\cyrfam=\titlecyr
    \scriptfont\cyrfam=\twelvecyr
    \scriptscriptfont\cyrfam=\tencyr%%
  \setbox\strutbox=\hbox{\vrule
      height12.3pt depth5.54pt width0pt}%
  \baselineskip=16pt\rm}

 %%%% RUNNING HEADINGS
\newbox\AuthorBox\newbox\TitleBox
\newbox\TFLinebox
\newbox\FLinebox
\newbox\HLinebox
\def\SetTFLinebox#1{\setbox\TFLinebox=\hbox{#1}}
\def\SetFLinebox#1{\setbox\FLinebox=\hbox{#1}}
\def\SetHLinebox#1{\setbox\HLinebox=\hbox{#1}}

 \def\SetAuthorHead#1{%
     \setbox\AuthorBox=\hbox{\ninepoint \it 
           \ignorespaces\frenchspacing#1\unskip}}
 \def\SetTitleHead#1{%
     \setbox\TitleBox=\hbox{\ninepoint \it
           \ignorespaces\frenchspacing#1\unskip}}

 %% Italic Spacing Correction
  \def\itSpacing{\relax}
  \def\itSpacingOff{\relax}

  %% Main section headings

 \def\Hrule{\hrule width0pt height0pt}

 %% skip used around proclamations, after section headings,
  % and before subsection-headings:
  \newskip\ProcSkip \ProcSkip 8pt plus2pt minus2pt

 \newskip\LastSkip
 \def\SaveLastSkip{\LastSkip\lastskip}
 \def\RestoreLastSkip{\vskip-\LastSkip\vskip\LastSkip}

 %% Do not indent next paragraph after a header:
 \def\NoindentAfter{\everypar={\setbox0=\lastbox\everypar={}}}

 \long\def\H#1\par#2\par{\notenumber=0 \titlepagetrue%
    {
    \baselineskip=20pt
    \parindent=0pt\parskip=0pt\frenchspacing
    \leftskip=0pt plus .2\hsize minus .3\hsize
    \rightskip=0pt plus .2\hsize minus .3\hsize
 \def\\{\unskip\break}%
    \pretolerance=10000 \Hfont #1\unskip\break
     \vskip7pt\Hrule
\hfill \Authorfont #2\hfill\hfill\unskip}
    \vskip48pt plus 4pt minus 4pt% 60pt=48+12pt
    \par\NoindentAfter\rm}

 \long\def\Hi#1\par#2\par{\notenumber=0 \titlepagetrue%
    {  \baselineskip=0pt  \parindent=0pt\parskip=0pt\frenchspacing
    \leftskip=0pt plus .2\hsize minus .3\hsize
    \rightskip=0pt plus .2\hsize minus .3\hsize
}
    \rm}

 %%% Minor section headings

 \newdimen\PageRemainder
  \def\SetPageRemainder{%\maxdimen case at page tops 12-91 LS
     \PageRemainder=\pagegoal
     \ifdim\PageRemainder=\maxdimen\PageRemainder=\vsize
     \else\advance\PageRemainder by -1\pagetotal\fi}

  \def\Rpt@{}\def\Rpt@@{}

  \long\def\HH#1\par{\par%A
  \SaveLastSkip\removelastskip\goodbreak
  \ifdim\LastSkip<30pt %24pt
     \LastSkip 30pt%24pt 
plus 3pt minus 2pt\fi
  \SetPageRemainder\advance\PageRemainder-\LastSkip
  \ifdim\PageRemainder<150pt
       \edef\Rpt@{remain = \the\PageRemainder\noexpand\\
                pagetotal=\the\pagetotal\noexpand\\
                           pagegoal=\the\pagegoal}%
          \fi
   \ifdim\PageRemainder<65pt %%Head plus 4 lines (approx)
       \ifdim\PageRemainder > 0pt
          \edef\Rpt@@{\noexpand\\
                      Had HH PageRemainder$<$\relax 65pt\noexpand\\
                      Hence forced break!}%
     \vskip 0pt plus .2\PageRemainder\eject %% Pull it out a bit
    \fi\fi
    \vskip\LastSkip\Hrule %%%%%%%%\Hrule added
    \pretolerance=10000\rightskip=0pt plus 3em%B
    \hangafter1 \hangindent=2.2em%
    \noindent
    \HHfont \unskip \Ednote{\Rpt@\Rpt@@}%
            \def\Rpt@{}\def\Rpt@@{}%
            \ignorespaces
            #1\par\rightskip=0pt\pretolerance=\StdPretolerance%
    \NoindentAfter
\tenpoint\rm%
     \medskip \vskip\ProcSkip}%interlineskip adds 2pt to this

  \long\def\HHH#1\par{\par%
  \SaveLastSkip\removelastskip\goodbreak
  \ifdim\LastSkip<\ProcSkip%
     \LastSkip\ProcSkip\fi
  \SetPageRemainder\advance\PageRemainder-\LastSkip
  \ifdim\PageRemainder<150pt
       \edef\Rpt@{remain = \the\PageRemainder\noexpand\\
                pagetotal=\the\pagetotal\noexpand\\
                           pagegoal=\the\pagegoal}%
       \fi
   \ifdim\PageRemainder<48pt  %% 4 lines
        \ifdim\PageRemainder > 0pt
             \edef\Rpt@@{\noexpand\\
                      Had HHH PageRemainder$<$\relax48pt\noexpand\\
                      Hence forced break!}%
       \vskip 0pt plus .2\PageRemainder\eject %% Pull it out a bit
      \fi\fi
   \vskip\LastSkip\par\noindent
   \HHHfont \unskip\Ednote{\Rpt@\Rpt@@}%
  \def\Rpt@{}\def\Rpt@@{}%
  \ignorespaces
   #1\unskip.\quad\rm\ignorespaces
   \ignorepars}

  \long\def\ignorepars#1\par{\def\Test{#1}%
     \ifx\Test\Empty\def\This{\ignorepars}%
        \else\def\This{\Test\par}\fi
           \This}
  \def\Empty{}

 \def\Abstract#1\par{\bgroup\Smallfonts\narrower\HHH #1\par}
 \def\endAbstract{\par\egroup}

 %%%%% Proclamations %%%%%

 \def\ProcBreak{\par%
    \ifdim\lastskip<8pt%
    \removelastskip%
    \penalty-200\vskip\ProcSkip\fi}

 \def\th#1\par{\ProcBreak \noindent
   {\thfont\ignorespaces
    #1\unskip.}\it\itSpacing\kern.4em\ignorepars}%\everymath{\/}

 \def\endth{\ProcBreak\rm\itSpacingOff }%\everymath{}

  %% the theorem statement will be in italic by default

 \def\pf#1\par{\ProcBreak %
    \noindent\pffont#1\unskip.\rm\itSpacingOff{\kern .7em}\ignorepars}

  %% \qed is alternative

  %% A Box for the QED
  \def\qedbox{\hbox{\vbox{
    \hrule width0.2cm height0.2pt
    \hbox to 0.2cm{\vrule height 0.2cm width 0.2pt
             \hfil\vrule height0.2cm width 0.2pt}
    \hrule width0.2cm height 0.2pt}\kern1pt}}

  %% Typing in \qed makes the qedbox right justified:
  \def\qed{\ifmmode\qedbox
    \else\unskip\ \hglue0mm\hfill\qedbox\ProcBreak\fi}

  \def \rk #1\par{\ProcBreak
     \noindent{\rkfont\ignorespaces #1\unskip.}%
     \rm\kern.6em\ignorepars}

  \def \df #1\par{\ProcBreak
     \noindent{\dffont\unskip\ignorespaces #1\unskip.}%
     \rm\kern.6em\ignorepars}

  \def \enddf {\medskip\ProcBreak }

  \def \eg #1\par{\ProcBreak
     \noindent\egfont\unskip\ignorespaces #1\unskip.
     \rm\kern.6em\ignorepars}

  \newdimen\Overhang

   \def\MaxTag@#1#2#3#4#5{\setbox0=\hbox{#4\ignorespaces#2\unskip}%
     \dimen0=\wd0\advance\dimen0 by#3
     \ifdim\dimen0<#5\relax\dimen0=#5\fi
     \expandafter\edef\csname #1Hang\endcsname{\the\dimen0}}

 \def\MaxItemTag#1{\MaxTag@{Item}{#1}{.4em}{\ItemStyle}{\parindent}}%
 \def\MaxItemItemTag#1{%
        \MaxTag@{ItemItem}{#1}{.4em}{\ItemItemStyle}{\parindent}}
 \def\MaxNrTag#1{\MaxTag@{Nr}{#1}{.5em}{\NrStyle}{\parindent}}
 \def\MaxReferenceTag#1{%
        \MaxTag@{Reference}{[#1]}{.6em}{\ninerm}{\parindent}}
 \def\MaxFootTag#1{\MaxTag@{Foot}{#1}{.4em}{\ninerm}{\z@}}

  %% \SetOverhang@ will prevent for tag-text collision
  \def\SetOverhang@{\Overhang=.8\dimen0%
     \advance\Overhang by \wd0\relax%nec!
     \ifdim\Overhang>\hangindent\relax%nec!
       \advance\Overhang by .25\dimen0%
       \Ednote{Tag is pushing text.}\osumess{Tag is pushing text.}%
     \else\Overhang=\hangindent
     \fi}

   %%% \Item
   \def\Item#1{\par\noindent
      \hangafter1\hangindent=\ItemHang
      \setbox0=\hbox{\ItemStyle\ignorespaces#1\unskip}%
      \dimen0=.4em\SetOverhang@% dimen0 is extra space
      \rlap{\box0}\kern\Overhang\ignorespaces}

   %%% \ItemItem
   \def\ItemItem#1{\par\noindent
      \hangafter1\hangindent=\ItemItemHang
      \setbox0=\hbox{\ItemItemStyle\ignorespaces#1\unskip}%
      \dimen0=.4em\SetOverhang@
      \advance\hangindent by \ItemHang
      \kern\ItemHang\rlap{\box0}%
      \kern\Overhang\ignorespaces}

  %%%% \Nr Items without hanging indentation
  \def\Nr#1{\par\noindent\hangindent=\NrHang %not really a hang
    \setbox0=\hbox{\NrStyle\ignorespaces#1\unskip}%
    \dimen0=.5em\SetOverhang@% dimen0 is extra space
    \rlap{\box0}\kern\Overhang
    \hangindent=\z@\ignorespaces}

  %%%% Roster (not compulsory)
  %%  endRoster has to remember \lastskip (e.g. from a \qed) through \egroup.
   \newskip\Rosterskip\Rosterskip 1pt plus1pt %% modifiable
   \def\Roster{\par\ifdim\lastskip<\Rosterskip\removelastskip\vskip\Rosterskip\fi
    \bgroup}
   \def\endRoster{\par\global\edef\LastSkip@{\the\lastskip}\removelastskip
       \egroup\penalty-50\LastSkip\LastSkip@\relax
       \ifdim\LastSkip<\Rosterskip\LastSkip\Rosterskip\fi
       \vskip\LastSkip}%%changed Feb/5/92 WN

 %%%%% Emphasis %%%%%

 %%%%% Vertical spacing %%%%%

 %%%%% References %%%%%

 \def\cite#1{%\relaxnext@
    \def\nextiii@##1,##2\end@{{\frenchspacing\rm 
      \lBr\ignorespaces##1\unskip{\rm,~\ignorespaces##2}\rBr}}%
    \IN@0,@#1@%
    \ifIN@\def\next{\nextiii@#1\end@}\else
    \def\next{{\rm\lBr#1\rBr}}\fi\next}

 %%%%% Bibliography %%%%%

   \def \Bib#1\par{%
       \par\removelastskip\SetPageRemainder
       \ifdim\PageRemainder < 97pt
        \ifdim\PageRemainder > 0pt
        \vfill\eject
       \fi\fi
    \ProcBreak \par\begingroup\parskip=0 pt%
    \goodbreak \vskip 15 pt plus 10 pt
    \noindent\null\hfill\Bibfont% \kern??pt%  (center over what?)
      \ignorespaces #1\unskip\hfill\null\par 
    \frenchspacing \Smallfonts\rm
    \parskip=2.5 pt plus 1 pt minus.5pt%
    \nobreak\vskip 12pt plus 2pt minus2pt\nobreak
    \leftskip=0 pt \baselineskip=10.5pt}

 \def\ReferenceTagSlide{0em}
  \def\ReferenceTagGap{.5em}

  \def \rf#1{\par\noindent
     \hangafter1\hangindent=\ReferenceHang      
     \setbox0=\hbox{\ninerm[\ignorespaces#1\unskip]}%        
     \dimen0=\ReferenceTagGap\SetOverhang@
     \rlap{\kern\ReferenceTagSlide\box0}%       
     \kern\Overhang\ignorespaces}

  \def\ref#1\par#2\par#3\par#4\par{%
     \rf{#1}#2\unskip,\ #3\unskip,\
     #4\unskip.}

  \def\endBib{\par\endgroup\vskip 12pt minus 6pt }

 %%%%% Coordinates %%%%%

  \long\def\Coordinates#1\endCoordinates{%\relax}
 {\par\vskip4pt\def\\{\unskip, }\Coordfont\baselineskip10.5pt\noindent#1}}

 \def\pagecontents{%\TRMargIns new, \Pagetot@l
  \gdef\Pagetot@l{\pagetotal}
  \ifvoid\TRMargIns\else
    \rlap{\kern\hsize\kern10pt\vbox to 0pt{%
         \box\TRMargIns\vss}}\fi
  \ifvoid\topins\else\unvbox\topins\fi
   \dimen@=\dp\@cclv \unvbox\@cclv % open up \box255
   \ifvoid\footins\else % footnote info is present
     \vskip\skip\footins
     \footnoterule
     \unvbox\footins\fi
   \ifr@ggedbottom \kern-\dimen@ \vfil \fi}

  %%%%% Some math accents %%%%%

 \newcount\Ht %pg121; Height register, used in Linefigure & accents

 \def \Acc{\expandafter } %%% What is this for?? WN

 \def\swthat{\raise -1.1 ex\hbox{\sam$\widehat{}$}}
 \def\swttilde{\raise -1.2 ex\hbox{\sam$\widetilde{}$}}
 \def \overdot{{\raise .2 ex \hbox to 0pt {\hss\bf\smash{.}\hss}}}
 \def \overcircle{{\raise .1 ex \hbox to 0pt
    {\sam$\eightpoint\scriptstyle\hss\circ\hss$}}}

 \def \Mathaccent#1#2{{\sam % E.g. #1=\widehat
  \setbox4=\hbox{$\vphantom{#2}$}
  \Ht=\ht4 %pg120
  \setbox5=\hbox{${#1}$}
  \setbox6=\hbox{${#2}$}
  \setbox7=\hbox to .5\wd6{}
  \copy7\kern .1\Ht \raise\Ht sp\hbox{\copy5}\kern-.1\Ht
  \copy7\llap{\box6}
  }}

  \def\SwtCheck #1{
        \ifmmode \check{#1}%
                \else \v {#1}%
                \fi}

 %%  \barpartial : bar over partial is common, tailor!
 \def\barpartial {%
   \kern .17 em
    \overline {\kern -.17 em\partial\kern-.03 em}%
    \kern .03 em}

 %%%   BEtter overline
 
  \def\Overline#1{\setbox1=\hbox{\sam ${#1}$}%
      \ifdim \wd1 > 6pt
    \kern .11 em
    \overline {\kern -.11 em#1\kern-.14 em}
    \kern .14 em
  \else
    \kern .03 em
    \overline {\kern -.03 em#1\kern-.04 em}
    \kern .04 em
  \fi}

 \def\SOverline#1{\setbox1=\hbox{\sam ${#1}$}%
      \ifdim \wd1 > 7pt
    \kern .22 em
    \overline {\kern -.22 em#1\kern-.09 em}%
    \kern .09 em
  \else
    \kern .10 em
    \overline {\kern -.10 em#1\kern-.04 em}%
    \kern .04 em
  \fi}

  %%% Better underline

 \def\Underline#1{\setbox1=\hbox{\sam ${#1}$}%
      \ifdim \wd1 > 6pt
    \kern .11 em
    \underline {\kern -.11 em#1\kern-.14 em}
    \kern .14 em
  \else
    \kern .03 em
    \underline {\kern -.03 em#1\kern-.04 em}
    \kern .04 em
  \fi}

 \def\SUnderline#1{\setbox1=\hbox{\sam ${#1}$}%
      \ifdim \wd1 > 7pt
    \kern .04 em
    \underline {\kern -.04 em#1\kern-.2 em}%
    \kern .2 em
  \else
    \kern .0 em
    \underline {\kern -.0 em#1\kern-.15 em}%
    \kern .15 em
  \fi}

  %%%%% Miscellaneous %%%%%

 \def \Blackbox
   {\leavevmode\hskip .3pt \vbox
   {\hrule height 5pt\hbox{\hskip 4.5pt}}\hskip .5pt}

 \def \XX{\Blackbox\kern.5pt\Blackbox} %% editorial use

  \def\.{.\kern1pt}

  %% unbreakable hyphen (by local change of hyphenchar to -1)
    \def\Hyphen{\edef\this{\the\hyphenchar\font}%
          \hyphenchar\font=-1\char\this\hyphenchar\font=\this}

  %% Prose In Math or Display 
 \ifx\undefined\text
  \def\text#1{\hbox{\rm #1}}\fi %% AMSTeX is more sophisticated

  %% Math Object Names (multi-character math object names)
  %%\nolimits can be cancelled
                                     % by a following \limits if wanted

%%%% Larry's mathsurround hacks:

   \everymath{}  %% initially, but later ...

  \def\PassMath@@{\aftergroup\AfterMath@} %% use \aftergroup LS 5-92

 \let\PassMath@\PassMath@@

 \def\AfterMath@{\futurelet\next\AfterMathMole@}

 \def\AfterMathMole@{%\show\next
      \ifcat\next\space% picks off CR and \par cases too; not \dots
          \def\this{}%{(space)}%
      \else
      \ifcat\next\egroup %
        \def\this{\osumess{Handset mathsurround?? ---(see dollar brace)}}%
      \else
      \def\this{\AAfterMath@}% this minority case slow
      \fi\fi
      \this}

 \def\hyphen@{-}
 \def\paren@{)}
 \def\apostr@{'}

 \def\MSC#1{\kern-.8\mathsurround#1\kern.8\mathsurround}

 \def\AAfterMath@#1{\def\Next{#1}%\show\Next%
    \IN@0\Next @,.;:!?\relax @%
    \ifIN@\def\this{\MSC{\Next}}%
    \else
    \ifx\Next\hyphen@\def\this{\futurelet\next\AfterHyphen@}%
    \else
    \ifx\Next\paren@\def\this{#1}%
    \else 
    \ifx\Next\apostr@\def\this{#1}%
    \else \def\this{\osumess{Handset mathsurround??}%
                 #1}\fi\fi\fi\fi
    \this}

 \def\AfterHyphen@#1{\def\Next{#1}%
   \ifx\Next\hyphen@\def\this{--}\else
   \ifcat\next\space%
   \def\this{\kern-\mathsurround\kern.05em- \Next}\else
   \def\this{\kern-\mathsurround\kern.05em\Hyphen\Next}\fi\fi\this}

%%%% switches
 \def\sam{\mathsurround=\z@\let\PassMath@\relax}  %
 \def\mas{\mathsurround=\StdMathsurround\let\PassMath@\PassMath@@}
 
 \def\Mas{\mathsurround=\StdMathsurround
                \everymath{\PassMath@}\let\PassMath@\PassMath@@}

 \def\m@th{\mathsurround=\z@\everymath{}}%% good general measure

 \def\m@@th{\mathsurround=\z@\everymath={}\let\m@th\relax}

\def\underbar#1{$\setbox\z@\hbox{#1}\dp\z@\z@
      \m@th \underline{\box\z@}$\relax}

\def\mathhexbox#1#2#3{\leavevmode
  \hbox{\m@@th$\m@th \mathchar"#1#2#3$}}

\def\dots{\relax\ifmmode\ldots\else$\m@th\ldots\,$\relax\fi}
   %%% this first \relax is ONLY original

\def\dotfill{\cleaders\hbox{\m@@th$\m@th \mkern1.5mu.\mkern1.5mu$}\hfill}
\def\rightarrowfill{$\m@th\mathord-\mkern-6mu%
  \cleaders\hbox{\m@@th$\mkern-2mu\mathord-\mkern-2mu$}\hfill
  \mkern-6mu\mathord\rightarrow$\relax}
\def\leftarrowfill{$\m@th\mathord\leftarrow\mkern-6mu%
  \cleaders\hbox{\m@@th$\mkern-2mu\mathord-\mkern-2mu$}\hfill
  \mkern-6mu\mathord-$\relax}

\def\downbracefill{$\m@th\braceld\leaders\vrule\hfill\braceru
  \bracelu\leaders\vrule\hfill\bracerd$\relax}
\def\upbracefill{$\m@th\bracelu\leaders\vrule\hfill\bracerd
  \braceld\leaders\vrule\hfill\braceru$\relax}

\def\angle{{\vbox{\m@@th\ialign{$\m@th\scriptstyle##$\crcr
      \not\mathrel{\mkern14mu}\crcr
      \noalign{\nointerlineskip}
      \mkern2.5mu\leaders\hrule height.34pt\hfill\mkern2.5mu\crcr}}}}

\def\big#1{{\m@@th\hbox{$\left#1\vbox to8.5\p@{}\right.\n@space$}}}
\def\Big#1{{\m@@th\hbox{$\left#1\vbox to11.5\p@{}\right.\n@space$}}}
\def\bigg#1{{\m@@th\hbox{$\left#1\vbox to14.5\p@{}\right.\n@space$}}}
\def\Bigg#1{{\m@@th\hbox{$\left#1\vbox to17.5\p@{}\right.\n@space$}}}
\def\n@space{\nulldelimiterspace\z@ \m@th}

\def\root#1\of{\setbox\rootbox\hbox{\m@@th$\m@th\scriptscriptstyle{#1}$}
  \mathpalette\r@@t}
\def\r@@t#1#2{\setbox\z@\hbox{\m@@th$\m@th#1\sqrt{#2}$\relax}
  \dimen@\ht\z@ \advance\dimen@-\dp\z@
  \mkern5mu\raise.6\dimen@\copy\rootbox \mkern-10mu \box\z@}

\def\mathph@nt#1#2{\setbox\z@\hbox{\m@@th$\m@th#1{#2}$}\finph@nt}

\def\mathsm@sh#1#2{\setbox\z@\hbox{\m@@th$\m@th#1{#2}$}\finsm@sh}

\def\@vereq#1#2{\lower.5\p@\vbox{\m@@th\baselineskip\z@skip\lineskip-.5\p@
    \ialign{$\m@th#1\hfil##\hfil$\crcr#2\crcr=\crcr}}}

\def\mathpalette#1#2{\sam\mathchoice{#1\displaystyle{#2}}%
  {#1\textstyle{#2}}{#1\scriptstyle{#2}}{#1\scriptscriptstyle{#2}}\mas}

\def\widehat#1{\setbox\z@\hbox{\sam$#1$}%
 \ifdim\wd\z@>\tw@ em\mathaccent"0\msbfam@5B{#1}%
 \else\mathaccent"0362{#1}\fi}
\def\widetilde#1{\setbox\z@\hbox{\sam$#1$}%
 \ifdim\wd\z@>\tw@ em\mathaccent"0\msbfam@5D{#1}%
 \else\mathaccent"0365{#1}\fi}

 \def\dots{\relax{}
  \ifmmode\def\thedots{\mdots@}\else\def\thedots{\tdots@}\fi %
  \thedots}

 %% \eqno and \leqno need protection
 \let\@ldeqno\eqno\let\@ldleqno\leqno
 \def\eqno{\everymath{}\@ldeqno} \def\leqno{\everymath{}\@ldleqno}

  \let\@ldeqalignno\eqalignno
  \def\eqalignno#1{\sam\@ldeqalignno{#1}\mas}
  \let\@ldeqalign\eqalign
  \def\eqalign#1{\sam\@ldeqalign{#1}\mas}

 \def\overrightarrow#1{\vbox{\m@th\ialign{##\crcr
      \rightarrowfill\crcr\noalign{\kern-\p@\nointerlineskip}
      $\hfil\displaystyle{#1}\hfil$\crcr}}}
 \def\overleftarrow#1{\vbox{\m@th\ialign{##\crcr
      \leftarrowfill\crcr\noalign{\kern-\p@\nointerlineskip}
      $\hfil\displaystyle{#1}\hfil$\crcr}}}
 \def\overbrace#1{\mathop{\vbox{\m@th\ialign{##\crcr\noalign{\kern3\p@}
      \downbracefill\crcr\noalign{\kern3\p@\nointerlineskip}
      $\hfil\displaystyle{#1}\hfil$\crcr}}}\limits}
 \def\underbrace#1{\mathop{\vtop{\m@th\ialign{##\crcr
      $\hfil\displaystyle{#1}\hfil$\crcr\noalign{\kern3\p@\nointerlineskip}
      \upbracefill\crcr\noalign{\kern3\p@}}}}\limits}

  \let\@ldmatrix\matrix
  \let\end@ldmatrix\endmatrix
  \def\matrix{\sam\@ldmatrix}
  \def\endmatrix{\end@ldmatrix\mas}
  \let\@ldgather\gather
  \let\end@ldgather\endgather
  \def\gather{\sam\@ldgather}
  \def\endgather{\end@ldgather\mas}
  \let\@ldalign\align
  \let\end@ldalign\endalign
  \def\align{\sam\@ldalign}
  \def\endalign{\end@ldalign\mas}
  \let\@ldaligned\aligned
  \let\end@ldaligned\endaligned
  \def\aligned{\sam\@ldaligned}
  \def\endaligned{\end@ldaligned\mas}
  \let\@ldtag\tag
  \def\tag{\sam\@ldtag}
   %
  %%% Commutative diagrams : use LamsCD too?

   \let\MinCDArrowWidth\minCDaw@

  %% will be redefined by BoxedEPS.tex

  %%%%% \FigureTitle %%%%%

%%%% End of Larry's mathsurround stuff
%%%% Start of Walter's insert corrections

\newskip\insertskipamount\newskip\inserthardskipamount
\insertskipamount 6pt plus2pt %This is medskipamount without shrink
\inserthardskipamount 6pt
\def\insertskip{\vskip\insertskipamount}
\newcount\SplitTest%        will be set to -1 if a topinsert has split
\def\SetSplitTest{\SplitTest\insertpenalties
  \insert\topins{\floatingpenalty1}%
  \advance\SplitTest-\insertpenalties}
\def\midinsert{\par
 \SaveLastSkip\penalty-150\SetSplitTest\RestoreLastSkip
 \ifnum\SplitTest=-1
  \@midfalse\p@gefalse\else\@midtrue\fi\@ins}
\def\@ins{\par\begingroup\setbox\z@\vbox\bgroup%
  \vglue\inserthardskipamount}
\def\endinsert{\egroup % finish the \vbox
  \if@mid \dimen@\ht\z@ \advance\dimen@\dp\z@
    \advance\dimen@\insertskipamount%            was 12pt (wn)
    \advance\dimen@\pagetotal\advance\dimen@-\pageshrink
    \ifdim\dimen@>\pagegoal\@midfalse\p@gefalse\fi\fi
  \if@mid%
    \ifdim\lastskip<\insertskipamount\removelastskip\insertskip\fi
    \nointerlineskip\box\z@\penalty-200\insertskip
  \else%
    \SaveLastSkip%                                  added (wn)
    \insert\topins{\penalty100 % floating insertion
    \splittopskip\z@skip
    \splitmaxdepth\maxdimen \floatingpenalty\z@
    \ifp@ge \dimen@\dp\z@
    \vbox to\vsize{\unvbox\z@\kern-\dimen@}% depth is zero
    \else \box\z@\nobreak\insertskip\fi}% was \bigskip\fi (wn)
    \RestoreLastSkip%                               added (wn)
   \fi\endgroup}
%% End Walter's insert stuff

 %%%%% Footnotes %%%%%

  \newcount\notenumber
  
  \def\note{\advance\notenumber by 1
    \footnote{\the\notenumber)}}

  \newbox\footbox

 %% The following modifies Plain TeX definitions, qv
  \def\footnote#1{\let\@sf\empty
    %{(the text)} is read later
    \ifhmode\edef\@sf{\spacefactor\the\spacefactor}\/\fi
    \sam${}^{\fam0 #1}$\@sf\vfootnote{#1}}%

  \def\vfootnote#1{\insert\footins\bgroup
     \interlinepenalty100 \splittopskip=1pt
     \floatingpenalty=20000
     \leftskip=0pt\rightskip=0pt%
     \parindent=.3em%% adjust
     \Smallfonts\rm%%osudeG added \Smallfonts
     \FootItem@{#1}%\strut% not nec
     \futurelet\next\fo@t}

  \def\FootItem@#1{\par\hangafter1\hangindent=\FootHang
     \setbox0=\hbox{\ignorespaces#1\unskip}%
     \dimen0=.4em\SetOverhang@% dimen0 is extra space
     \noindent\rlap{\box0}\kern\Overhang\ignorespaces}

  %\MaxFootTag{2)}%% in param file

  \def\fo@t{\ifcat\bgroup\noexpand\next \let\next\f@@t
    \else\let\next\f@t\fi \next}
  \def\f@@t{\bgroup\aftergroup\@foot\let\next}
  \def\f@t#1{\baselineskip=10pt\lineskip=1pt
            \lineskiplimit=0pt #1\@foot}%
     %%osudeG added \baselineskip=? pt\lineskiplimit=0pt
  \def\@foot{%%% special strut osu for end of each note
        \hbox{\vrule height0pt depth5pt width0pt}
        \egroup}
  \skip\footins=12 pt plus 0pt minus 0pt %% was \bigskipamount
    %% space added when footnote is present
  \count\footins=1000 % footnote magnification factor (1 to 1)
  \dimen\footins=8in % maximum footnotes per page

 %%%% Altenatives

  %%  Editorial stuff (delete??)

 \def\osumess#1{\EdSpider{\immediate\write16{Line \the\inputlineno: #1}}}%
 \def\HideEdStuff{\gdef\EdSpider##1{}}

 \font\BigSym=cmmi10 scaled \magstep 4

 \def\change{\InLMargin{\hbox{\BigSym \char63\kern10pt}}}

 \def\beginchange{\InLMargin{\hbox{\sam\twelvepoint$\heartsuit$\kern10pt}}}

 \def\endchange{\InLMargin{\hbox{\sam\twelvepoint$\spadesuit$\kern10pt}}}

 \def\InLMargin#1{\strut\vadjust{%
     \kern-\strutdepth
     \vtop to \strutdepth{%
         \baselineskip\strutdepth
         \llap{\sam$\smash{\hbox{\EdSpider{#1}}}$}\null}}}

 \def\strutdepth{\dp\strutbox}
 \def\strutheight{\ht\strutbox}

 \def\NoteInRMargin#1{\strut\vadjust{%
     \kern-1.001\strutdepth
     \vtop to \strutdepth{%
       \baselineskip\strutdepth
       \vss\rlap{\ninepoint\unskip\hskip\hsize
         \vtop to 0pt{%
           \hsize=16em\hfuzz=\hsize
           \leftskip=10pt%
           \rightskip=0pt plus 10000pt%
           \baselineskip=9.8pt\lineskip=.2pt%
           \let\\\break
           \noindent\EdSpider{#1}\vss}%
                \kern10pt}\hbox{}}%%\hbox{}=\null crucial!!
       }}

 \def\ednote#1{\NoteInRMargin{\tentt #1}}

 \def\cbar{\InLMargin{%
      \dimen0=\strutdepth\advance\dimen0 by \lineskip
      \vrule width 3pt
      height \strutheight depth \dimen0 \kern
      3pt}}

 \def\ccbar{\InLMargin{%
      \dimen0=2\strutdepth\advance\dimen0 by 2\lineskip
      \vrule width 3pt
        height 3\strutheight depth \dimen0 \kern
      3pt}}

 \newinsert\TRMargIns
 \dimen\TRMargIns=\maxdimen
 %\count\TRMargIns=0
 %\skip\TRMargIns=0pt

  \def\Ednote#1{\insert\TRMargIns{%
       \vbox to 0pt{\hsize=140pt\hfuzz=\hsize
           \leftskip=6pt%
           \rightskip=0pt plus 10000pt%
           \baselineskip=9.8pt\lineskip=.2pt%
           \let\\\break
           %\vglue\pagetotal% misplaces notes if inserts are present
           \SetPageRemainder% This ...
           \vglue540pt\vglue-\PageRemainder%  .. is a fix (WN)
           \noindent\EdSpider{\tentt #1}\vss}%
       \smallskip}}

 \def\KillEdStuff{\def\ednote##1{}\def\Ednote##1{}%
      \let\change\relax\let\beginchange\relax\let\endchange\relax
       \let\cbar\relax\let\ccbar\relax}

 %%% Compatibility with osumrip.sty
  %%

 %%% Parameters
  \topskip=12pt
  \newskip\StdBaselineskip % to set \baselineskip
  \StdBaselineskip 12pt
  \lineskip=1.1pt
  \lineskiplimit=.8pt
  \widowpenalty=10000 % 8000 to 10000
  \clubpenalty=10000  % 8000 to 10000
  \abovedisplayskip=6pt plus 1pt minus 1pt
  \abovedisplayshortskip=3pt plus 1.5pt
  \belowdisplayskip=6pt plus 1pt minus 1pt
  \belowdisplayshortskip=5pt plus 1pt minus 1pt
  \hfuzz=1.5pt   % Enable overfull box warnings at console

  \def\StdPretolerance{100}
  \tolerance=\StdPretolerance

  \newdimen\StdMathsurround
  \StdMathsurround=1.5pt % 1pt usual without \Mas
  \mathsurround=\StdMathsurround
  \Mas                   %% sophisticated mathsurround on
 % \Sam                   %% sophisticated mathsurround off

%% marker before English punctuation in displayed math
   \def\prose{\relax\hbox{\kern.6\StdMathsurround}}
  
  \def\StdParskip{0pt}    %% Larry wants {2pt plus 1pt}
  \parskip=\StdParskip
  \parindent=0.5cm
 
%%%% load Times for main body font

  \def\Times{ptmr  } 
  \def\TimesI{ptmri  } 
  \def\TimesB{ptmb  }
  \def\TimesBI{ptmbi  }
  \def\HelveticaN{phvrrn }

  =\Times at 10bp% roman text
  =\TimesB at 10bp% boldface extended
   % slanted roman
  \font\tenit=\TimesI at 10bp% text italic
  =\TimesBI at 10bp

  \font\tenmrm=cmr10  %%new name for math role at full size

%%%%% Fonts at ninepoint %%%%%

    =\Times at 9bp 
    \font\nineit=\TimesI at 9bp 
    =\TimesB at 9bp 
    =\TimesBI at 9bp 

    =\HelveticaN at 9bp 
       % see below

%%%%% Fonts at twelvepoint %%%%%

  =\Times at 12bp
  \font\twelveit=\TimesI at 12bp
  =\TimesB at 12bp

%%%%% Fonts at titlepoint %%%%%

  \font\titleit=\TimesI at 14.4bp
  =\TimesB at 14.4bp

 \SetAuthorHead{AuthorHead} % needs \ninepoint since box set
 \SetTitleHead{TitleHead}  % notably \HeaderFont

%%%% Char adjustments %%%%

  \def\lBr{\raise.125ex\hbox{[\kern.1125ex}}
  \def\rBr{\raise.125ex\hbox{\kern.1125ex]}}

 \setbox\footbox=\hbox{\Smallfonts 2)~}

%% Some optional font dimension and spacing 
%% adjustments beyond this point

%% Correct the lousy spacing of italic f (a hack).

  \bgroup
  \catcode`\@=11 %localised
  \gdef\itSpacing{%
     \xspaceskip=.31em plus.1em minus.05em \sfcode `f=2001
     \itWarning@\let\itWarning@\itWarning@@}
  \gdef\itSpacingOff{%
     \xspaceskip=0pt \sfcode `f=1000
     \let\itWarning@\relax}
   \global\let\itWarning@\relax
  \gdef\itWarning@@{\errmessage{%
  Special italic spacing already in force
  (you have probably omitted an ``endth'').
  See itSpacing macro in osuPSfnt.sty
         }}
  \egroup

 %%% Provisional fontdimen settings
  %%
 \fontdimen1\titlebf=0.0pt
 \fontdimen2\titlebf=3.6135pt
 \fontdimen3\titlebf=2.8908pt
 \fontdimen4\titlebf=1.44539pt
 \fontdimen5\titlebf=6.64882pt
 \fontdimen6\titlebf=14.45398pt
 \fontdimen7\titlebf=1.60439pt

 \fontdimen1\tenbi=0.26794pt
 \fontdimen2\tenbi=2.50937pt
 \fontdimen3\tenbi=2.00749pt
 \fontdimen4\tenbi=1.00374pt
 \fontdimen5\tenbi=4.59717pt
 \fontdimen6\tenbi=10.03749pt
 \fontdimen7\tenbi=1.11415pt

 \fontdimen1\twelverm=0.0pt
 \fontdimen2\twelverm=3.01125pt
 \fontdimen3\twelverm=2.409pt
 \fontdimen4\twelverm=1.2045pt
 \fontdimen5\twelverm=5.39615pt
 \fontdimen6\twelverm=12.045pt
 \fontdimen7\twelverm=1.33699pt

 \fontdimen1\twelveit=0.27731pt
 \fontdimen2\twelveit=3.01125pt
 \fontdimen3\twelveit=2.409pt
 \fontdimen4\twelveit=1.2045pt
 \fontdimen5\twelveit=5.37207pt
 \fontdimen6\twelveit=12.045pt
 \fontdimen7\twelveit=1.33699pt

 \fontdimen1\twelvebf=0.0pt
 \fontdimen2\twelvebf=3.01125pt
 \fontdimen3\twelvebf=2.409pt
 \fontdimen4\twelvebf=1.2045pt
 \fontdimen5\twelvebf=5.5407pt
 \fontdimen6\twelvebf=12.045pt
 \fontdimen7\twelvebf=1.33699pt

 \fontdimen1\tenrm=0.0pt
 \fontdimen2\tenrm=2.50937pt
 \fontdimen3\tenrm=2.00749pt
 \fontdimen4\tenrm=1.00374pt
 \fontdimen5\tenrm=4.49678pt
 \fontdimen6\tenrm=10.03749pt
 \fontdimen7\tenrm=1.11415pt

 \fontdimen1\tenit=0.27731pt
 \fontdimen2\tenit=2.50937pt
 \fontdimen3\tenit=2.00749pt
 \fontdimen4\tenit=1.00374pt
 \fontdimen5\tenit=4.47672pt
 \fontdimen6\tenit=10.03749pt
 \fontdimen7\tenit=1.11415pt

 \fontdimen1\tenbf=0.0pt
 \fontdimen2\tenbf=2.50937pt
 \fontdimen3\tenbf=2.00749pt
 \fontdimen4\tenbf=1.00374pt
 \fontdimen5\tenbf=4.61723pt
 \fontdimen6\tenbf=10.03749pt
 \fontdimen7\tenbf=1.11415pt

 \fontdimen1\ninerm=0.0pt
 \fontdimen2\ninerm=2.25842pt
 \fontdimen3\ninerm=1.80673pt
 \fontdimen4\ninerm=0.90337pt
 \fontdimen5\ninerm=4.0471pt
 \fontdimen6\ninerm=9.03374pt
 \fontdimen7\ninerm=1.00273pt

 \fontdimen1\nineit=0.27731pt
 \fontdimen2\nineit=2.25842pt
 \fontdimen3\nineit=1.80673pt
 \fontdimen4\nineit=0.90337pt
 \fontdimen5\nineit=4.02904pt
 \fontdimen6\nineit=9.03374pt
 \fontdimen7\nineit=1.00273pt

 \fontdimen1\ninebf=0.0pt
 \fontdimen2\ninebf=2.25842pt
 \fontdimen3\ninebf=1.80673pt
 \fontdimen4\ninebf=0.90337pt
 \fontdimen5\ninebf=4.15552pt
 \fontdimen6\ninebf=9.03374pt
 \fontdimen7\ninebf=1.00273pt

 %%% \SetExtraSpaces \MaxSpaceFactor \SetSpaceFactors
  %%  See TeXbook, page 76.

 \newcount\MaxSpaceFactor
 \MaxSpaceFactor=3000 %% to reset later

 %%%%% Tag styles and (hang-) indents
 \def\ItemStyle{\rm}
 \def\NrStyle{\rm}
 \def\ItemItemStyle{\rm}

 %% Analog dimensioning, convenient for local modifications:
 \MaxItemTag{(iii)}
 \MaxItemItemTag{(iii)}
 \MaxNrTag{(2)}
 \MaxFootTag{2)}
 % \MaxReferenceTag{AaaAA} % for biblio
 \def\ReferenceHang{30pt}

 \catcode`\@=\active

%%%%% End of hack of Neumann-Siebenmann macros

\loadbold

=\Times  
=\Times scaled750
=\Times scaled650
\font\rms=\Times scaled 920 

=\TimesBI scaled 860
=\TimesI scaled 860

\textfont0=\rrm  
\scriptfont0=\erm 
\scriptscriptfont0=\srm

\def\Augment#1#2{%
    \toks0\expandafter{#1}\toks2{#2}%
    \edef#1{\the\toks0\the\toks2}}

 \font\twelverma=\Times  scaled 1200
 \font\tenrma=\Times  scaled 1000
 \font\ninerma=\Times scaled 920
 =\Times scaled 840
 \font\sevenrma=\Times scaled 760
 =\Times scaled 680
 \font\fiverma=\Times scaled 600

 \Augment\tenpoint{%
  \textfont0=\tenrma  \scriptfont0=\sevenrma  
  \scriptscriptfont0=\fiverma  }

 \Augment\ninepoint{%
  \textfont0=\ninerma  \scriptfont0=\sevenrma 
  \scriptscriptfont0=\fiverma}

 \Augment\twelvepoint{%
  \textfont0=\twelverma  \scriptfont0=\ninerma  
  \scriptscriptfont0=\sevenrma}

\mathsurround=1pt
\hsize=13.45truecm
\vsize=19.5truecm
\hoffset=1.25truecm
\voffset=2truecm
\advance\baselineskip by 2pt

\predefine\til{\~}
\def\~#1{\relax\ifmmode\widetilde{#1}\else\til{#1}\fi}

\redefine \le{\leqslant}

\define \wt#1{\mathaccent"0365{#1}}
\define \wh#1{\mathaccent"0362{#1}}

\define \iss{\,\Mathaccent{\raise -.8 ex\hbox{$\widetilde{}$\kern.1em}}\rightarrow\,}

\define \pr{\operatorname{\fam0 pr}}

\Mas
\HideEdStuff
\rm 
 
%%%% For GT headers and footers:

\def\issn{{\nineit ISSN 1464-8997 (on line) 1464-8989 (printed)}}

\def\gtp{{\nineit Published 10 December 2000: \ \copyright\ Geometry \& 
Topology Publications}}

\def\gtv3{{\nineit Geometry \& Topology Monographs, Volume 3 (2000) --
Invitation to higher local fields}}

%%%%% For section idents:

\def\lione
{{\rms Geometry \& Topology Monographs}}

\def \litwo{{\rms Volume 3: Invitation to higher local fields
}} 

\def\tinfo #1.#2.#3-#4
{{
\noindent  {\lione} \hfill 
\par 
\vskip-1.5pt
\noindent {\litwo} \hfill
\par 
\vskip-1,5pt
\noindent {\rms Part #1, section #2, pages #3--#4} \hfill
\vskip24pt 
}}

\def\tinfos #1.#2.#3-#4
{{
\noindent  {\lione} \hfill 
\par 
\vskip-1.5pt
\noindent {\litwo} \hfill
\par 
\vskip-1.5pt
\noindent {\rms Pages #3--#4} \hfill
\vskip24pt 
}}

\def\tinfoi #1
{{
\noindent  {\lione} \hfill 
\par 
\vskip-1.5pt
\noindent {\litwo} \hfill
\par 
\vskip-1.5pt
\noindent {\rms Pages iii--xi: Introduction and contents} \hfill
\vskip26pt 
}}

%%%% Set headers and footers %%%%

  \def\titlepagehead{\hfil}

  \newif\iftitlepage\titlepagefalse
  \newif\ifblankpage\blankpagefalse
  \def\makeheadline{
     \ifblankpage{}\else%
     \iftitlepage
\vbox{\line{\vbox to 8.5pt{}
\ninerm
\copy\HLinebox \hfill
\hglue5mm\ninebf\folio 
\titlepagehead}}%
      \else
\vbox{\ifodd\pageno\rightheadline\else\leftheadline\fi}%
      \fi\vskip 12pt\fi}%
     \def\rightheadline{\line{\vbox to 8.5pt{}%
      \ninerm
\copy\TitleBox \hfill
\hglue5mm\ninebf\folio}}%
     \def\leftheadline{\line{\vbox to 8.5pt{}%
        \unskip\ninerm\unskip\ninebf\folio\hglue5mm
      %*%
 \hfill \copy\AuthorBox
%\hfill
}}

 \footline={\ifblankpage{}\else
\iftitlepage\ninepoint\sam\hfill%} 
\line{\vbox to 8.5pt{}%\ninerm
\copy\TFLinebox
\hfill
\hglue5mm %\ninebf\folio
}
            \else
\ninepoint\sam\hfill%}
\line{\vbox to 8.5pt{}%\ninerm
\copy\FLinebox
\hfill 
\hglue5mm
}
\hfil\fi\global\titlepagefalse\fi}

\def\blankpage{{\blankpagetrue\noindent\hbox to 10pt{\hss}\vfill
\pagebreak}}

\tenpoint\rm %% always start here
 
  %%% all done and macros loaded!

\input xy
\xyoption{all}

\input picture
\makeatletter
\def\@inmatherr#1{\relax}
\resetatcatcode
\def\makebox(0,0)[cc]#1{\vbox to 0pt{\vss \hbox to 0pt{\hss #1\hss}\vss}}

\pageno=223

\tinfo II.3.223-237

\SetTFLinebox{\gtp }
\SetFLinebox{\gtv3 }
\SetHLinebox{\issn}

\H 3. The Bruhat--Tits buildings over \\ 
higher dimensional local fields  

A. N. Parshin

\SetAuthorHead{A. N. Parshin}
\SetTitleHead{Part II. Section 3. 
The Bruhat--Tits buildings over 
higher dimensional local fields \qquad\qquad}

\HH 3.0. Introduction

A  generalization  of  the  Bruhat--Tits
     buildings for the groups $PGL(V)$ over $n$-di\-me\-n\-sional  local  fields
was introduced in \cite{P1}.
The
main object of the  classical Bruhat--Tits theory  is  a  simplicial  complex  attached  to  any
reductive  algebraic  group  $G$  defined over a field $K$. There are two
parallel theories in  the case where  $K$  has no   additional  structure
or  $K$ is a local (or more generally, complete discrete valuation) field.
  They are known as the {\it spherical} and  {\it  euclidean}
buildings  correspondingly  (see subsection~3.2 for
a brief introduction,  \cite{BT1}, \cite{BT2} for original papers and
\cite{R}, \cite{T1} for the surveys).

In the generalized theory of buildings they correspond to local fields of dimension  zero and of 
dimension one. 
 The construction of the Bruhat--Tits  building  for  the  group $PGL(2)$ 
over two-dimensional  local field was described in detail in \cite{P2}.  
Later 
V.~Ginzburg and M. Kapranov   extended the theory to 
arbitrary reductive groups over a two-dimensional local fields \cite{GK}. Their
definition coincides with ours for $PGL(2)$ and is different for higher ranks.
But it seems that they are closely related (in  the case of the groups of
type $A_l$). It remains to develop the theory for arbitrary reductive groups over
local fields of dimension greater than two. 

In this work we describe the structure of the higher building for the
group $PGL(3)$ over a two-dimensional local field. 
We refer to \cite{P1}, \cite{P2} for the motivation of these constructions.
    
This work contains four subsections. In~3.1 
we collect  facts about   the Weyl group. Then in~3.2 we briefly describe the building 
for $PGL(2)$ over a local field of dimension not greater than two;
for details see \cite{P1}, \cite{P2}.
In~3.3 we study 
the building 
for $PGL(3)$ over a local field $F$ of dimension one  
and in~3.4 we 
describe  the building over a two-dimensional local field. 

%These notes were partially written up 
% during  my visit to the Abdus Salam ICTP in the autumn of 1997. 
%I am very much grateful to M. S. Narasimhan
%for the excellent working conditions. 
%My special thanks go to J. Juyumaya who kindly produced the diagrams. 

\smallskip

We use the notations of section~1 of Part~I. 

If $K$ is an $n$-dimensional local field,  
let $\Gamma_K$ be the valuation group of the
discrete valuation of rank $n$ on $K^*$;
the choice of a system of local parameters $t_1,\dots,t_n$ of $K$
induces an isomorphism of $\Gamma_K$ and
the lexicographically ordered group $\Bbb Z^{\oplus n}$.

Let $K$ ($K=K_2$, $K_1$, $K_0=k$) be a two-dimensional local field.  
Let $O=O_K$, $M=M_K$, $\Cal O=\Cal O_K$, $\Cal M=\Cal M_K$ (see subsection~1.1
of Part~I). 
Then $O=\pr^{-1}(\Cal O_{K_1})$, $M=\pr^{-1}(\Cal M_{K_1})$
where $\pr\colon \Cal O_K\to K_1$ is the residue map.
Let $t_1,t_2$ be a system of local parameters of $K$. 

If $K \supset {\Cal O}$ is the fraction field of a ring ${\Cal O}$ we call
${\Cal O}$-submodules $J \subset K$ fractional ${\Cal O}$-ideals (or
simply fractional ideals). 

The ring $O$  has the following
     properties: 
\Roster   
\Item{(i)} $ O/M  \simeq  k,\quad K^{*}  \simeq  \langle t_1\rangle\times
\langle t_2\rangle  \times O^*,\quad O^{*}  \simeq  k^{*}\times (1 +  M)$;

\Item{(ii)} every finitely  generated  fractional  $O$-ideal    is  
principal  and equal to  
$$P(i, j) = (t_1^{i}t_2^{j})\quad \text{ for some $i,j\in \Bbb Z$} $$
(for the notation $P(i, j)$ see loc.cit.); 

\Item{(iii)} every infinitely generated fractional ${O}$-ideal  is equal to
$$P(j) = \Cal M_K^j=\langle t_1^{i}t_2^{j}:i\in \Bbb Z\rangle\quad
\text{for some $j \in \Bbb Z$}$$ 
(see \cite{FP}, \cite{P2} or section 1 of Part~I). 
The set of these ideals is totally ordered with respect to 
the inclusion.
\endRoster 

\vskip 1cm

\HH 3.1. The Weyl group

Let $B$ be the image of 
$$  
\pmatrix            O & O& \dots & O\\
                    M & O& \dots & O\\
                      & & \dots &           \\
                    M & M     &  \dots & O
            \endpmatrix 
$$
in   $PGL(m, K)$.  
 Let  $N$  be  the
subgroup of monomial matrices.

\df Definition 1

                   Let 
$ T = B \bigcap N $
be the image of
$$\pmatrix                          
 O^{*} & \dots & 0 \\
                                               & \ddots &  \\
                                            0   & \dots & O^{*}
                                            \endpmatrix 
$$
in $G$.

The group 
$$W=W_{K/K_1/k} = N/T$$ 
is called the
{\it Weyl group}. 
\enddf

There is a rich structure of subgroups in $G$ which have  many common 
properties with the theory of BN-pairs. 
In particular, there are  
Bruhat, Cartan and Iwasawa decompositions (see \cite{P2}). 

     The Weyl group $W$  contains the following elements of order two
$$ s_{i} = \pmatrix
                   1  & \dots  &  0  &    &   & 0 & \dots & 0 \\
                      & \ddots &     &    &   &   &     &    \\
                   0  & \dots  &  1  &    &   &   &     &  0 \\
                   0  & \dots  &     & 0  & 1 &   & \dots  & 0     \\
                   0  &        &     & 1 & 0 &   &  \dots  & 0      \\
                   0  & \dots  &     &    &   & 1 & \dots   & 0  \\
                      &        &     &    &   &   & \ddots &       \\
                   0  &  \dots &  0  &    &   & 0  & \dots  & 1
                   \endpmatrix, \quad i = 1, ..., m - 1;
$$
$$ w_{1} = \pmatrix
                   0 & 0 & \dots & 0 & t_1 \\
                   0 & 1 & \dots & 0 & 0 \\
                     &   & \dots &   &    \\
                     &   & \dots &   &   \\
                     &   & \dots & 1 & 0   \\
         t_1^{-1} & 0 & \dots & 0 & 0
         \endpmatrix
,\quad w_{2} = \pmatrix
                   0 & 0 & \dots & 0 & t_2 \\
                   0 & 1 & \dots & 0 & 0 \\
                     &   & \dots &   &    \\
                     &   & \dots &   &   \\
                     &   & \dots & 1 & 0   \\
         t_2^{-1} & 0 & \dots & 0 & 0
         \endpmatrix. 
$$

The  group $W$ has the following properties:
\Roster
 \Item{(i)}    
      $W$ is generated by the set $S$ of its elements of order two,

  \Item{(ii)}  there is an exact sequence 
        $$
        0 \rightarrow E \rightarrow W_{K/K_1/k}
       \rightarrow  W_{K}  \rightarrow  1, 
       $$
        where
       $ E$ is the kernel of the addition map
$$\underbrace{\Gamma_K\oplus\dots \oplus \Gamma_K}_{\text{$m$ times}}\to \Gamma_K $$ 
      and $ W_{K}$ is isomorphic to the symmetric group
$S_{m}$;  
     
      \Item{(iii)}  the elements $s_{i}$, $i = 1, \dots, m - 1$ define a splitting of
the exact sequence and the subgroup
 $\langle s_{1}, \dots, s_{m - 1}\rangle $
acts on $E$ by permutations.
 \endRoster

In contrast with the situation in the theory of BN-pairs
the pair $(W,S)$ is  not a Coxeter group
and furthermore there is no subset $S$ of involutions in  $W$  such  that
$(W,S)$ is a Coxeter group (see \cite{P2}).
\vskip 1cm

\HH 3.2. Bruhat--Tits building for $PGL(2)$ over a local field 
\unskip\break \phantom{}\enspace of dimension $\le 2$

\phantom{}\par

\smallskip 

In this subsection we briefly recall the main constructions.
For more details see \cite{BT1}, \cite{BT2}, \cite{P1}, \cite{P2}. 

\HHH 3.2.1

Let $k$ be a field (which can be viewed as a 0-dimensional local field).
Let $V$ be a vector space over $k$ of dimension two.

\df Definition 2

The {\it spherical building} of $PGL(2)$ over $k$
is a zero-dimensional complex $$\Delta(k)=\Delta(PGL(V),k)$$ whose vertices
are  lines in $V$.

The group $PGL(2,k)$ acts on $\Delta(k)$ transitively.
The Weyl group (in this case it is of order two)
acts on $\Delta(k)$ and its orbits are {\it apartments}
of the building.
\enddf

\HHH 3.2.2

Let $F$ be a complete discrete valuation field with residue field $k$.
Let $V$ be a vector space over $F$ of dimension two.
We say that $L \subset V$ is a lattice if $L$ is an ${\Cal O}_F$-module.
Two submodules $L$ and $L'$ belong to the same class $\langle L\rangle \sim 
\langle L'\rangle $ if and only if $L = aL'$, with $a \in F^{*}$.

\df Definition 3

The {\it euclidean building} of $PGL(2)$ over $F$
is a one-dimensional complex $\Delta(F/k)$ whose vertices 
are equivalence classes $\langle L \rangle$
of lattices.
Two classes $\langle L \rangle$
and $\langle L' \rangle$
are connected by an edge if and only if
for some choice of $L, L'$ there is an exact sequence
$$0\to L'\to L\to k\to 0.$$
Denote by $\Delta_i(F/k)$ the set of $i$-dimensional simplices
of the building $\Delta(F/k)$.
\enddf

The following {\it link property} is important:
\Roster
\Item{} Let $P\in \Delta_0(F/k)$ be represented by a lattice $L$.
Then the link of $P$ ($=$ the set of edges of $\Delta(F/k)$
going from $P$) is in one-to-one correspondence
with the set of  lines in the vector space
$V_P=L/\Cal M_F L$ (which is
$\Delta(PGL(V_P),k)$).
\endRoster

 The orbits of the Weyl group $W$ (which is in this case an infinite group with two generators of order two)  are infinite sets consisting of $x_i=\langle L_i \rangle$,
$L_i=\Cal O_F\oplus \Cal M_F^i$.

$$
\picture(300,10)

\put(120,0){\circle*{5}}
\put(150,0){\circle*{5}}
\put(180,0){\circle*{5}}

\put(110,0){\line(1,0){80}}

\put(98,0){\makebox(0,0)[cc]{$\dots$}}
\put(202,0){\makebox(0,0)[cc]{$\dots$}}

\put(120,8){\makebox(0,0)[cc]{$x_{i-1}$}}
\put(150,10){\makebox(0,0)[cc]{$x_{i}$}}
\put(180,10){\makebox(0,0)[cc]{$x_{i+1}$}}

\endpicture
$$

\medskip

An element $w$ of the Weyl group acts in the following way:
if $w\in E=\Bbb Z$ then $w$ acts by translation of even length;
if $w\not\in E$ then $w$ acts as an involution with a unique fixed point
$x_{i_0}$: 
$w(x_{i+i_0})=x_{i_0-i}$.

To formalize the connection of $\Delta(F/k)$ with $\Delta(F)$
we define a {\it boundary} point of $\Delta(F/k)$ as a class of
half-lines such that the intersection of every two half-lines
from the class is a half-line in both of them.
The set of the boundary points is called the {\it boundary} of
$\Delta(F/k)$.

There is an isomorphism between $PGL(2,F)$-sets
$\Delta(F)$ and the boundary of $\Delta(F/k)$:
if a half-line is represented by $L_i=\Cal O_F\oplus \Cal M_F^i$, $i>0$,
then the corresponding vertex of $\Delta(F)$ is the line
$F\oplus (0)$ in $V$.

It seems reasonable to slightly change the notations to make
the latter isomorphisms more transparent.

\df Definition 4 (\cite{P1})

Put $\Delta_{.} [0](F/k)=$ the complex of classes of $\Cal O_F$-submodules
in $V$ isomorphic to $F\oplus \Cal O_F$ (so $\Delta_{.} [0](F/k)$ is isomorphic to 
$\Delta(F)$)
 and put
$$\Delta_{.} [1](F/k)=\Delta(F/k).$$
Define the {\it building of $PGL(2)$ over $F$}
as the union
$$\Delta_{.}(F/k)=\Delta_{.}[1](F/k)\bigcup \Delta_{.}[0](F/k)$$
and call the subcomplex $\Delta_{.} [0](F/k)$
the {\it boundary} of the building.
The discrete topology on the boundary can be extended to the whole building.
\enddf

\HHH 3.2.3

Let $K$ be a two-dimensional local field.

Let $V$ be a vector space over $K$ of dimension two.
We say that $L \subset V$ is a lattice if $L$ is an ${O}$-module.
Two submodules $L$ and $L'$ belong to the same class $\langle L\rangle \sim 
\langle L'\rangle $ if and only if $L = aL'$, with $a \in K^{*}$.

\df Definition 5 (\cite{P1})

Define the vertices of the building
of $PGL(2)$ over $K$
as
$$
\aligned
\Delta_0[2](K/K_1/k)&=\text{classes of $O$-submodules $L\subset V$: $L\simeq O\oplus O$} 
\\
\Delta_0[1](K/K_1/k)&=\text{classes of $O$-submodules $L\subset V$: $L\simeq O\oplus \Cal O$} 
\\
\Delta_0[0](K/K_1/k)&=\text{classes of $O$-submodules $L\subset V$: $L\simeq O\oplus K$.} 
\endaligned
$$
Put
$$\Delta_0(K/K_1/k)= \Delta_0[2](K/K_1/k)\bigcup \Delta_0[1](K/K_1/k)\bigcup \Delta_0[0](K/K_1/k).$$

A set of $\{L_\alpha\}$, $\alpha\in I$,  of $O$-submodules in $V$ is called 
a {\it chain} if 
\Roster
\Item{(i)} for every $\alpha \in I$ and for every $a \in K^{*}$ there exists an 
$\alpha' \in I$ such that $aL_{\alpha} = L_{\alpha'}$,

\Item{(ii)} the set $\{ L_{\alpha}, \alpha \in I \}$ is totally ordered by the
inclusion.
\endRoster 
A chain $\{ L_{\alpha}, \alpha \in I \}$ is  called a {\it maximal chain} if
 it cannot
be included in a strictly larger set satisfying the same conditions (i) and
(ii). 

We say that $\langle L_{0}\rangle , \langle L_{1}\rangle, \dots , \langle L_{m}
\rangle $ belong to a {\it simplex} of
dimension 
$m$ if and only if  the $L_{i}$, $i = 0, 1, ..., m$ belong to a maximal chain of
${\Cal O}_F$-submodules in $V$. The faces and the degeneracies can be defined in
a standard way (as a deletion or repetition of a vertex). 
See \cite{BT2}.
\enddf

 Let $\{L_{\alpha} \}$ be a maximal chain of $O$-submodules in
the space $V$. There are exactly three types of maximal chains (\cite{P2}):
\Roster   
\Item{(i)} if the chain contains a module $L$ isomorphic to $O \oplus
O $ then all the modules of the chain are
of that type and the chain is uniquely determined by its segment 
$$\dots \supset O \oplus O  \supset
 M \oplus O \supset
 M \oplus M  \supset \dots .$$
 
   \Item{(ii)} if the chain contains a module $L$ isomorphic to $O \oplus \Cal O$ then  
 the chain can be restored from 
the segment: 
$$\dots\supset O\oplus \Cal O\supset O\oplus P(1,0)\supset O\oplus P(2,0)\supset\dots\supset O\oplus\Cal M\supset\dots$$
(recall that $P(1,0)=M$).
    
    \Item{(iii)}  if the chain contains a module $L$ isomorphic to $\Cal O \oplus \Cal O$ then  
 the chain can be restored from 
the segment: 
$$\dots\supset \Cal O\oplus\Cal O \supset P(1,0)\oplus \Cal O\supset P(2,0)\oplus \Cal O\supset\dots\supset \Cal M\oplus \Cal O\supset\dots .$$
 \endRoster

\vskip 1cm

\HH 3.3. Bruhat--Tits building for $PGL(3)$ over a local field $F$ 
\unskip\break \phantom{}\enspace of dimension 1

Let $G = PGL(3)$.

Let $F$ be a one-dimensional local field, $F \supset {\Cal O}_F \supset \Cal M_F$, 
${\Cal O}_F/\Cal M_F \simeq k$.
%Put $O_F=\Cal O_F$, M_F=\Cal M_F$. 

Let  $V$ be a vector space over $F$ of dimension three. 
Define lattices in $V$ and their equivalence similarly to the definition
of 3.2.2. 
%We say that $L \subset V$ is a lattice if $L$ is an ${\Cal O}_F$-module.
%Two submodules $L$ and $L'$ belong to the same class $\langle L\rangle \sim 
%\langle L'\rangle $ if and only if $L = aL'$, with $a \in F^{*}$. 

First we define the vertices of the building and then the simplices. The
result will be a simplicial set $\Delta_{.}(G, F/k)$. 

\df Definition 6

 The {\it  vertices} of the Bruhat--Tits building:
  $$
\aligned
 \Delta_{0}[1](G, F/k) &= \{ \text{classes of $\Cal O_F$
-submodules $ L \subset
V$ :  $L \simeq  {\Cal O}_F \oplus {\Cal O}_F \oplus {\Cal O}_F $}\}, \\
\Delta_{0}[0](G, F/k) & = \{ \text{classes of $\Cal O_F$
-submodules $ L \subset V$
: 
 $L \simeq {\Cal O}_F \oplus {\Cal O}_F \oplus F$}\\
&\quad \text{  or
$L \simeq {\Cal O}_F \oplus F \oplus F$}  \}, \\
\Delta_{0}(G, F/k) &= \Delta_{0}[1](G, F/k) \cup 
\Delta_{0}[0](G, F/k).
\endaligned
$$
We say that the points of $\Delta_{0}[1]$ are  {\it  inner}  points, 
the  points  of $\Delta_{0}[0]$ are  {\it  boundary} points.
 Sometimes we delete $G$ and $F/k$ from the notation  if  this
does not lead to  confusion.
\enddf

We  have defined the vertices only. For the simplices of higher dimension
we have the following: 

\df Definition 7

Let $\{ L_{\alpha}, \alpha \in I \}$ be a  set of ${\Cal O}_F$-submodules in
$V$. We say that $\{ L_{\alpha}, \alpha \in I \}$ is a {\it chain} if
\Roster
\Item{(i)} for every $\alpha \in I$ and for every $a \in K^{*}$ there exists an 
$\alpha' \in I$ such that $aL_{\alpha} = L_{\alpha'}$,

\Item{(ii)} the set $\{ L_{\alpha}, \alpha \in I \}$ is totally ordered by the
inclusion.
\endRoster 
A chain $\{ L_{\alpha}, \alpha \in I \}$ is  called a {\it maximal chain} if
 it cannot
be included in a strictly larger set satisfying the same conditions (i) and
(ii). 

We say that $\langle L_{0}\rangle , \langle L_{1}\rangle, \dots , \langle L_{m}
\rangle $ belong to a {\it simplex} of
dimension 
$m$ if and only if  the $L_{i}$, $i = 0, 1, ..., m$ belong to a maximal chain of
${\Cal O}_F$-submodules in $V$. The faces and the degeneracies can be defined in
a standard way (as a deletion or repetition of a vertex). 
See \cite{BT2}.
\enddf

To describe the structure of the building we first need  to determine all
types of the maximal chains. Proceeding as  in \cite{P2} (for 
$PGL(2)$) we get the following result.

\smallskip

\th Proposition 1

 There are exactly three types of   maximal chains of ${\Cal O}_F$-submodules in the space 
$V$:
\Roster
\Item{(i)} the chain contains a module isomorphic to ${\Cal O}_F \oplus
{\Cal O}_F \oplus {\Cal O}_F$. Then all the modules from the chain are
of that type and the chain has the following structure:
$$ \dots \supset \Cal M_F^i L \supset \Cal M_F^i L' \supset \Cal M_F^iL''
\supset \Cal M_F^{i + 1}L \supset \Cal M_F^{i + 1}L' \supset \Cal M_F^{i + 1} L''
\supset \dots $$ 
where $\langle L\rangle, \langle L'\rangle, \langle L''\rangle \in \Delta_{0}(G, F/k)[1] $ and 
     $L \simeq  {\Cal O}_F \oplus {\Cal O}_F \oplus {\Cal O}_F$, 

\Item{} $L' \simeq  {\Cal O}_F \oplus {\Cal O}_F \oplus \Cal M_F$, $L'' \simeq  {\Cal O}_F \oplus \Cal M_F
 \oplus \Cal M_F$. 
     
\Item{(ii)} the chain contains a module   isomorphic to ${\Cal O}_F \oplus {\Cal O}_F
\oplus F$. Then the chain has the following structure: 
$$~\dots \supset \Cal M_F^{i}L \supset \Cal M_F^{i}L' \supset \Cal M_F^{i +1}L 
     \supset \dots $$
      where $\langle L\rangle, \langle L'\rangle \in \Delta_{0}(G, F/k)[0] $ and
     $L \simeq  {\Cal O}_F \oplus {\Cal O}_F \oplus F$, 
$L' \simeq  \Cal M_F \oplus {\Cal O}_F \oplus F$. 
     
\Item{(iii)} the chain contains a  module  isomorphic to  ${\Cal O}_F \oplus F \oplus F$. Then 
the chain has the following structure: 
$$\dots \supset \Cal M_F^i L \supset \Cal M_F^{i+1}L \supset \dots$$ 
          where $\langle L\rangle\in \Delta_{0}(G, F/k)[0] $. 
   \endRoster
\endth

We see that the chains of the first type correspond to two-simplices, of the
second type --- to edges and the last type represent some vertices. It means that
the simplicial set
$\Delta_{.}$ is a  disconnected
union of its subsets $\Delta_{.}[m], ~m = 0, 1$. The dimension of the
subset $\Delta_{.}[m]$ is equal to one for $m = 0$ and to two for $m = 1$.

Usually  the  buildings  are  defined as combinatorial complexes having a
system of subcomplexes called apartments (see, for example,
\cite{R}, \cite{T1}, \cite{T2}). 
We show how to introduce them for the higher building.  

\df Definition 8

Fix a basis
$e_{1},  e_{2}, e_{3} \in V$.  The {\it apartment} defined by
this basis is the following set
$$ \Sigma_{.} = \Sigma_{.}[1] \cup \Sigma_{.}[0], $$
where
$$
\aligned
 \Sigma_{0}[1] = \{\langle L\rangle : L =& a_{1}e_{1} \oplus a_{2}e_{2} \oplus a_{3}e_{3}, \\
&\text{ where 
$a_{1}, a_{2}, a_{3}$ are  ${\Cal O}_F$-submodules in  $F$
isomorphic to $\Cal O_F$} \}\\
 \Sigma_{0}[0] = \{ \langle L\rangle : L =& a_{1}e_{1} \oplus a_{2}e_{2} \oplus a_{3}e_{3}, \\
&\text{where 
$a_{1}, a_{2}, a_{3}$ are $\Cal O_F$-submodules in $F$
 isomorphic
 either}\\
&\text{ to $\Cal O_F$ or to $F$}\\
&\text{and at least one $a_i$ is isomorphic to $F$}\}.
\endaligned
$$
$\Sigma_{.}[m] $ is the minimal subcomplex having $\Sigma_{0}[m]$ as vertices.
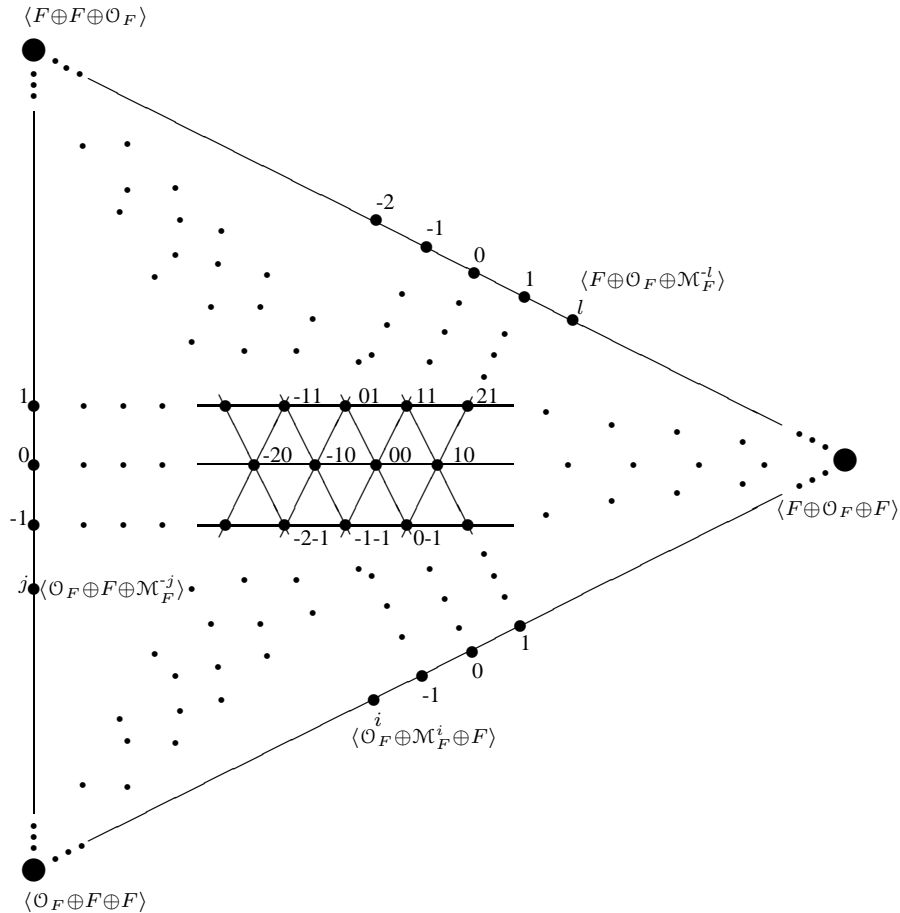
\enddf

\vglue-0.1in 
{\ninepoint
$$
\unitlength=0.75pt
\let\ss\scriptstyle    
\picture(400,425)\sam
\def\ccO{{\cal O}_F}
\def\ccM{{\cal M}_F}

\put(-5,53){\line(0,1){322}}
\put(20,40){\line(2,1){318}}
\put(20,390){\line(2,-1){318}}

\put(-5,27){\circle*{10}}
\put(-5,403){\circle*{10}}
\put(367,215){\circle*{10}}

\put(-5,37){\circle*{2}}
\put(-5,42){\circle*{2}}
\put(-5,47){\circle*{2}}
\put(-5,382){\circle*{2}}
\put(-5,387){\circle*{2}}
\put(-5,392){\circle*{2}}

\put(5,32){\circle*{2}}
\put(11,35){\circle*{2}}
\put(17,38){\circle*{2}}
\put(5,398){\circle*{2}}
\put(10,395){\circle*{2}}
\put(16,392){\circle*{2}}

\put(358,209){\circle*{2}}
\put(352,206){\circle*{2}}
\put(346,203){\circle*{2}}
\put(358,221){\circle*{2}}
\put(352,224){\circle*{2}}
\put(346,227){\circle*{2}}

\put(18,185){\circle*{2}}
\put(36,185){\circle*{2}}
\put(54,185){\circle*{2}}

\put(18,213){\circle*{2}}
\put(36,213){\circle*{2}}
\put(54,213){\circle*{2}}

\put(18,240){\circle*{2}}
\put(36,240){\circle*{2}}
\put(54,240){\circle*{2}}

\put(230,190){\circle*{2}}
\put(260,196){\circle*{2}}
\put(290,200){\circle*{2}}
\put(320,204){\circle*{2}}

\put(240,213){\circle*{2}}
\put(270,213){\circle*{2}}
\put(300,213){\circle*{2}}
\put(330,213){\circle*{2}}

\put(230,237){\circle*{2}}
\put(260,231){\circle*{2}}
\put(290,228){\circle*{2}}
\put(320,223.5){\circle*{2}}
\put(67.5,156){\circle*{2}}
\put(50.5,126){\circle*{2}}
\put(34.5,96){\circle*{2}}
\put(17.5,66){\circle*{2}}

\put(91.5,160){\circle*{2}}
\put(76.4,140){\circle*{2}}
\put(60,116){\circle*{2}}
\put(38,86){\circle*{2}}

\put(115.5,160){\circle*{2}}
\put(99,140){\circle*{2}}
\put(79.5,120){\circle*{2}}
\put(62,100){\circle*{2}}

\put(144,165){\circle*{2}}
\put(123,145){\circle*{2}}
\put(102,125){\circle*{2}}
\put(81,105){\circle*{2}}
\put(60,85){\circle*{2}}
\put(38,65){\circle*{2}}

\put(150,162){\circle*{2}}
\put(157,148){\circle*{2}}
\put(164,134){\circle*{2}}

\put(176.5,165){\circle*{2}}
\put(183.5,151){\circle*{2}}
\put(190,138){\circle*{2}}

\put(201.5,172){\circle*{2}}
\put(206,162){\circle*{2}}
\put(211,152){\circle*{2}}
\put(67.5,269){\circle*{2}}
\put(50.5,299){\circle*{2}}
\put(34.5,329){\circle*{2}}
\put(17.5,359){\circle*{2}}

\put(91.5,265){\circle*{2}}
\put(76.4,285){\circle*{2}}
\put(60,309){\circle*{2}}
\put(38,339){\circle*{2}}

\put(115.5,265){\circle*{2}}
\put(99,285){\circle*{2}}
\put(79.5,305){\circle*{2}}
\put(62,325){\circle*{2}}

\put(144,260){\circle*{2}}
\put(123,280){\circle*{2}}
\put(102,300){\circle*{2}}
\put(81,320){\circle*{2}}
\put(60,340){\circle*{2}}
\put(38,360){\circle*{2}}

\put(150,263){\circle*{2}}
\put(157,277){\circle*{2}}
\put(164,291){\circle*{2}}

\put(176.5,260){\circle*{2}}
\put(183.5,274){\circle*{2}}
\put(190,287){\circle*{2}}

\put(201.5,253){\circle*{2}}
\put(206,263){\circle*{2}}
\put(211,273){\circle*{2}}

\put(-10,10){%
{$\ss\langle \ccO\oplus F \oplus F \rangle$}}
\put(-10,416){%
{$\ss \langle F\oplus F \oplus \ccO \rangle$}}
\put(335,190){%
{$\ss \langle F\oplus \ccO\oplus F \rangle$}}

\put(-5,240){\circle*{5}}
\put(83,240){\circle*{5}}
\put(110,240){\circle*{5}}
\put(138,240){\circle*{5}}
\put(166,240){\circle*{5}}
\put(194,240){\circle*{5}}

\put(-12,241){$\ss1$}
\put(114,241){$\ss\text{-}11$}
\put(144,241){$\ss01$}
\put(170,241){$\ss11$}
\put(198,241){$\ss21$}

\put(-5,213){\circle*{5}}
\put(96,213){\circle*{5}}
\put(124,213){\circle*{5}}
\put(152,213){\circle*{5}}
\put(180,213){\circle*{5}}

\put(-12,214){$\ss0$}
\put(100,214){$\ss\text{-}20$}
\put(129,214){$\ss\text{-}10$}
\put(158,214){$\ss00$}
\put(187,214){$\ss10$}

\put(-5,185){\circle*{5}}
\put(83,185){\circle*{5}}
\put(110,185){\circle*{5}}
\put(138,185){\circle*{5}}
\put(166,185){\circle*{5}}
\put(194,185){\circle*{5}}

\put(-16,185){$\ss\text{-}1$}
\put(114,176){$\ss\text{-}2\text{-}1$}
\put(142,176){$\ss\text{-}1\text{-}1$}
\put(169,176){$\ss0\text{-}1$}

\put(-5,156){\circle*{5}}
\put(-12,156){$\ss j$}
\put(-3,153){%
{$\ss \langle \ccO\oplus F \oplus \ccM^{\text{-}j}\rangle$}}

\put(242,279){\circle*{5}}
\put(220,290){\circle*{5}}
\put(197,301){\circle*{5}}
\put(175,313){\circle*{5}}
\put(152,325){\circle*{5}}

\put(245,295){%
{$\ss \langle F\oplus \ccO \oplus \ccM^{\text{-}l}\rangle$}}
\put(244,282){$\ss l$}
\put(220,295){$\ss1$}
\put(197,306){$\ss0$}
\put(175,318){$\ss\text{-}1$}
\put(152,330){$\ss\text{-}2$}

\put(218,139){\circle*{5}}
\put(196,127){\circle*{5}}
\put(173,116){\circle*{5}}
\put(151,105){\circle*{5}}

\put(218,128){$\ss1$}
\put(196,115){$\ss0$}
\put(173,104){$\ss\text{-}1$}
\put(151,93){$\ss i$}
\put(140,85){%
{$\ss \langle \ccO\oplus \ccM^i \oplus F \rangle$}}

\put(70,240){\line(1,0){145}}
\put(70,213){\line(1,0){145}}
\put(70,185){\line(1,0){145}}

\put(80,180){\line(1,2){32}}
\put(108,180){\line(1,2){32}}
\put(136,180){\line(1,2){32}}
\put(164,180){\line(1,2){32}}

\put(80,245){\line(1,-2){32}}
\put(108,245){\line(1,-2){32}}
\put(136,245){\line(1,-2){32}}
\put(164,245){\line(1,-2){32}}

\endpicture
$$
}\vglue-0.1in

It can be shown that the building $ \Delta_{.}(G, F/k)$
is glued from the apartments, namely 

$$ \Delta_{.}(G, F/k) = \bigsqcup_{\text{all bases of $V$}} \Sigma_{.}\,/ \,
\text{an equivalence relation} $$
(see \cite{T2}).

We can  make  this description  more transparent by drawing all that in
the  picture above where the dots of different  kinds  belong  to  the
different parts of the building. In contrast with the case of the group
$PGL(2)$ it is not easy to draw the whole building and we restrict
ourselves to an apartment.

Here the inner vertices are represented by the lattices
$$ij = \langle {\Cal O}_F \oplus \Cal M_F^i \oplus \Cal M_F^j\rangle, \quad
i,j \in \Bbb Z.$$

     The definition of the boundary gives a topology on
$\Delta_{0}(G, F/k)$ which is discrete on  both  subsets  $\Delta_{0}[1]$
and $\Delta_{0}[0]$. The convergence of the inner points to the boundary
points is given by the following rules:
$$
\aligned
 \langle {\Cal O}_F \oplus \Cal M_F^i \oplus \Cal M_F^j\rangle 
& @>{j \rightarrow -\infty}>> 
\langle {\Cal O}_F \oplus \Cal M_F^i \oplus F\rangle , \\
\langle {\Cal O}_F \oplus \Cal M_F^i \oplus \Cal M_F^j\rangle 
& @>{j \rightarrow \infty}>> 
\langle F \oplus F \oplus {\Cal O}_F\rangle ,
\endaligned  $$
because $  \langle {\Cal O}_F \oplus \Cal M_F^i \oplus \Cal M_F^j\rangle = \langle \Cal M_F^{-j} \oplus \Cal M_F^{-j+i} \oplus
{\Cal O}_F\rangle$. The convergence in the other two directions can be defined
along the same line (and it is shown on the picture).
It is easy to
extend it to the higher simplices.

  Thus,   there is  the  structure of a simplicial 
topological space on the apartment and then we define it on the whole
building using the gluing procedure.
  This  topology  is  stronger than the
topology usually introduced to connect the inner part and the boundary  together. The connection with standard ``compactification'' of the building
is given by the following map:

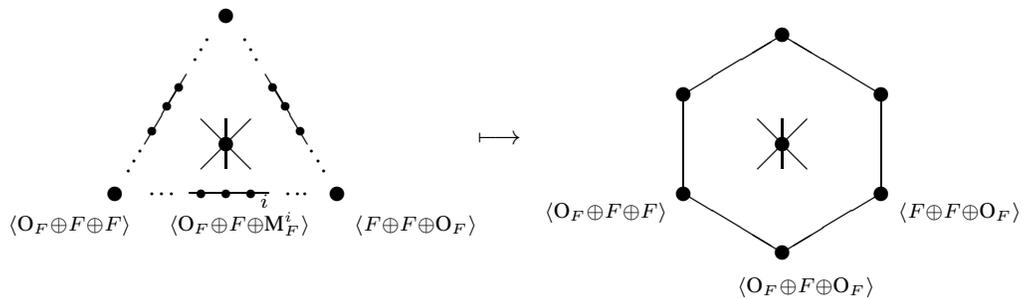
\vglue-0.4in
{\ninepoint
\def\ccO{\text{\cal O}_F}
\def\ccM{\text{\cal M}_F}

$$
\unitlength=0.85 pt
\let\ss\scriptstyle    
\picture(350,180)

%Hexagono
\put(250,120){\line(5,3){40}}
\put(330,80){\line(-5,-3){40}}

\put(330,80){\line(0,1){40}}
\put(250,80){\line(0,1){40}}

\put(330,120){\line(-5,3){40}}
\put(250,80){\line(5,-3){40}}

%
%cruz
\put(290,90){\line(0,1){20}}
\put(280,90){\line(1,1){20}}
\put(300,90){\line(-1,1){20}}

%puntos 
\put(250,120){\circle*{6}}
\put(330,80){\circle*{6}}
\put(330,120){\circle*{6}}
\put(250,80){\circle*{6}}
\put(290,56){\circle*{6}}
\put(290,100){\circle*{6}}
\put(290,144){\circle*{6}}
\put(192,70){%\scriptsize
{$\ss \langle {\ccO}\oplus F \oplus F \rangle$}}

\put(270,40){%\scriptsize
{$\ss \langle {\ccO}\oplus F \oplus {\ccO} \rangle$}}

\put(335,70){%\scriptsize
{$\ss \langle F\oplus F\oplus {\ccO} \rangle$}}

%fin

%triangulo

%\put(20,80){\line(1,0){87}}
%\put(20,80){\line(3,5){45}}
%\put(110,80){\line(-3,5){45}}

\put(50,80){\line(1,0){32}}
\put(32,100){\line(3,5){17}}
\put(98,100){\line(-3,5){17}}

%puntos grandes
\put(20,80){\circle*{6}}
\put(110,80){\circle*{6}}
\put(65,152){\circle*{6}}
\put(65,100){\circle*{6}}

%puntos medianos
\put(55,80){\circle*{3}}
\put(65,80){\circle*{3}}
\put(75,80){\circle*{3}}

\put(35,105){\circle*{3}}
\put(41,115){\circle*{3}}
\put(46,123){\circle*{3}}

\put(95,105){\circle*{3}}
\put(89,115){\circle*{3}}
\put(84,123){\circle*{3}}

%puntos pequegnos
\put(25.4,89){\circle*{1}}
\put(27.8,93){\circle*{1}}
\put(30.2,97){\circle*{1}}
\put(52.5,134){\circle*{1}}
\put(54.9,138){\circle*{1}}
\put(57.3,142){\circle*{1}}
\put(104.5,89){\circle*{1}}
\put(102.2,93){\circle*{1}}
\put(99.5,97){\circle*{1}}
\put(77.4,134){\circle*{1}}
\put(75,138){\circle*{1}}
\put(72.6,142){\circle*{1}}
\put(35,80){\circle*{1}}
\put(39,80){\circle*{1}}
\put(43,80){\circle*{1}}

\put(90,80){\circle*{1}}
\put(94,80){\circle*{1}}
\put(97,80){\circle*{1}}
%
%cruz
\put(65,90){\line(0,1){20}}
\put(55,90){\line(1,1){20}}
\put(75,90){\line(-1,1){20}}
\put(-25,65){%\scriptsize
{$\ss \langle {\ccO}\oplus F \oplus F \rangle$}}

\put(77,74){%\scriptsize
{$\ss i$}}

\put(40,65){%\scriptsize
{$\ss \langle {\ccO}\oplus F \oplus  \ccM^i \rangle$}}

\put(115,65){%\scriptsize
{$\ss \langle F\oplus F\oplus {\ccO} \rangle$}}
%fin
%flecha
\put(165,100){$\longmapsto$} 

\endpicture
$$}

\vglue-0.55in

This map is bijective on the inner simplices and on a part of the boundary 
can be described as
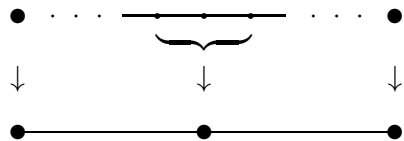
\vglue-0.3in
$$\unitlength0.8pt
\picture(325,80)
% ejes
%\put(0,0){\line(0,1){100}}
%\put(0,0){\line(1,0){400}}

%lineas

\put(125,60){\line(1,0){70}}

\put(80,60){\circle*{6}}

\put(95,60){\circle*{1}}
\put(105,60){\circle*{1}}
\put(115,60){\circle*{1}}

\put(140,60){\circle*{3}}
\put(160,60){\circle*{3}}
\put(180,60){\circle*{3}}

\put(207,60){\circle*{1}}
\put(217,60){\circle*{1}}
\put(227,60){\circle*{1}}

\put(242,60){\circle*{6}}

%parte de abajo
\put(80,10){\circle*{6}}
\put(160,10){\circle*{6}}
\put(242,10){\circle*{6}}

\put(80,10){\line(1,0){160}}
%
%flechas
\put(75,30){$\downarrow$}
\put(237,30){$\downarrow$}
\put(155,30){$\downarrow$}
%
%brace
\put(137.5,53){$\underbrace{\hbox to 1.15cm{\hss}}$}
\endpicture$$
\vglue-0.05in
\noindent
We note that the complex is not a CW-complex but only a closure finite 
complex. 
 This ``compactification'' was used by G. Mustafin \cite{M}.

We have two kinds of connections with the  buildings for other fields and
groups. First, for the local field $F$
there are two local fields of dimension 0, namely 
$F$ and $k$. Then 
for every  $P\in \Delta_{0}[1](PGL(V), F/k)$ the  $\text{Link}(P)$
is equal to  
$\Delta_{.}(PGL(V_{P}), k)$ 
where $V_P = L/\Cal M_FL$ if $P = \langle L \rangle$ and the $\text{Link} (P)$ is the 
boundary of the $\text{Star} (P)$. Since the apartments for the
$PGL(3, k)$ are  hexagons, we can also observe this property on the picture. 
The analogous relation with the building of $PGL(3, K)$ is more complicated.
It is shown on the picture above.

The other relations work if we change the group $G$ but not the field.
We see that three different lines go out from every inner point in
the apartment. They  represent the apartments of the group $PGL(2, F/k)$.
They correspond to different embeddings of the $PGL(2)$ into $PGL(3)$.

Also we can describe the action of the Weyl group $W$ on an apartment.
If we fix a basis, the extension
$$ 0 \rightarrow \Gamma_F \oplus \Gamma_F \rightarrow W \rightarrow
S_3 \rightarrow 1 $$ 
 splits. 
The elements from $S_3 \subset W$ act either as rotations around
the point $00$ or as reflections. The elements of  $\Bbb Z \oplus \Bbb Z \subset W $ can be represented as triples of integers (according to property (ii) in
the previous subsection). Then they correspond to translations of the lattice of
 inner points along the three directions going from the point $00$. 

If we fix an embedding $PGL(2) \subset PGL(3)$ then the apartments and
the Weyl groups are connected as follows:  
$$ 
\aligned 
 &\Sigma_{.}(PGL(2)) \subset \Sigma_{.}(PGL(3)), \\
&\CD
0 @>>>  \Bbb Z @>>> W' @>>> S_2 @>>> 1\\
@. @VVV @VVV @VVV @. \\
 0 @>>> \Bbb Z \oplus \Bbb Z @>>> W 
 @>>> S_3 @>>> 1
\endCD
\endaligned  
$$
where $W'$ is a Weyl group of the group $PGL(2)$ over the field $F/k$. 

\vskip 1cm

\HH 3.4. Bruhat--Tits building for $PGL(3)$ over a local field 
\unskip\break \phantom{}\enspace of  dimension 2

Let $K$ be a two-dimensional local field.  Denote by $V$ a vector
space over $K$ of dimension three.  Define lattices in V and their
equivalence in a similar way to 3.2.3.  We shall consider the 
following {\it types} of lattices: 
   $$  \matrix
   \Delta_0[2] & &  222 & \langle O \oplus O \oplus O \rangle  \\
      \Delta_0[1]  &  &   221   & 
\langle O \oplus O \oplus \Cal O \rangle  \\
& & 211 & \langle O \oplus \Cal O \oplus \Cal O \rangle
 \\
\Delta_0[0] & & 220  & \langle O \oplus  O \oplus K \rangle \\
 & & 200 & \langle O \oplus  K \oplus K \rangle 
\endmatrix
       $$ 
To define the buildings we repeat the procedure from the previous subsection.

\df Definition  10

 The {\it  vertices} of the Bruhat--Tits building 
     are the elements of the following set: 
  $$ \Delta_{0}(G, K/K_1/k) = \Delta_{0}[2] \cup \Delta_{0}[1]
\cup \Delta_{0}[0].$$
To define the simplices of higher dimension we can repeat word by word   
Definitions~7 and~8 of the previous subsection replacing  the ring ${\Cal O}_F$ by the ring $O$ (note that we work only with the types of lattices listed  above). We call the
     subset $\Delta[1]$ the {\it inner boundary} of the building and the
     subset $\Delta[0]$ the {\it external boundary}. The points in
      $\Delta[2]$ are the {\it inner} points.
  \enddf

To describe the structure of the building we first need to determine all
types of the maximal chains. Proceeding as  in \cite{P2} for
$PGL(2)$ we get the following result.

\th Proposition 2

 Let $\{L_{\alpha} \}$ be a maximal chain of $O$-submodules in
the space $V$. There are exactly five types of maximal chains:
 \Roster   
\Item{(i)} If the chain contains a module $L$ isomorphic to $O \oplus
O \oplus O$ then all the modules of the chain are
of that type and the chain is uniquely determined by its segment 
$$\dots \supset O \oplus O \oplus O \supset
 M \oplus O \oplus O \supset
 M \oplus M \oplus O \supset
 M \oplus M \oplus M \supset
\dots $$

   \Item{(ii)} If the chain contains a module $L$ isomorphic to $O \oplus O 
\oplus {\Cal O}$ then the chain can be restored from
the segment: 
\bigskip
\centerline{\let\ss\scriptstyle
\let\sss\scriptscriptstyle
\vbox{\hbox
{$\ss\dots\supset \text{``${\Cal O} \oplus {\Cal O} \oplus {\Cal O}$''} \supset
\dots \supset O \oplus O \oplus {\Cal O} \supset
M \oplus O \oplus {\Cal O} \supset
M \oplus M \oplus {\Cal O} \supset \dots \supset
\Cal M \oplus \Cal M\oplus {\Cal O}$}\vskip -6pt 
\hbox{$\sss\qquad\qquad\xymatrix{
*=0{}\ar@/_/  [rrrrrr]|{\sss\text{ quotient $\sss\simeq   K_1 \oplus  K_1$ }}  & & & & &  & *=0{}
}\quad$} 
\vskip 6pt
\hbox{$\ss=     %\underbrace
{\Cal M \oplus \Cal M \oplus {\Cal O} \supset
     \dots \supset \Cal M\oplus \Cal M\oplus O \supset \Cal M\oplus
     \Cal M\oplus M \supset \dots \supset \Cal M\oplus \Cal M
\oplus
\Cal M 
\supset \dots}
$}
\vskip -6pt 
 \hbox{$\sss\qquad\xymatrix{
*=0{}\ar@/_/  [rrrrr]|{\sss\text{ quotient $\sss\simeq   K_1$ }}  & & & & & *=0{}
}\qquad\quad$}}}
\bigskip 
\endRoster
 \Roster
\Item{}    Here the modules  isomorphic to  
${\Cal O} \oplus {\Cal O} \oplus {\Cal O}$
     do not belong to this chain  and are inserted as in  
the proof of
     Proposition~1 of \cite{P2}.
    
    \Item{(iii)}  All the modules $L_{\alpha} \simeq  O \oplus {\Cal O} \oplus
{\Cal  O}$. Then the chain contains a piece
\bigskip 
\centerline{\ninepoint
\def\Matt
{ \xymatrix{
*=0{}\ar@/_/  [rrrrrrrr]|{\text{ quotient $\simeq   K_1$ }}  &  & & & & & &  & *=0{}
} }
\def\Mattt
{ \xymatrix{
*=0{}\ar@/_/  [rrrrrr]|{\text{ quotient $\simeq   K_1$ }}  & & &  & &  & *=0{}
} }
\quad \vbox
{ \hbox
{$
\dots \supset  \text{``${\Cal O} \oplus {\Cal O} \oplus {\Cal O}$''} \supset
\dots \supset O \oplus {\Cal O} \oplus {\Cal O} \supset
M \oplus {\Cal O} \oplus {\Cal O} \supset \dots
 \supset
\Cal M
\oplus {\Cal O} \oplus {\Cal O}
$}
\vskip -6pt 
 \hbox{$\qquad\qquad\qquad\Matt
\qquad\quad$}  
\vskip 6pt 
\hbox
{$
 =  \Cal M\oplus {\Cal O} \oplus {\Cal O}
\supset 
     \dots \supset \Cal M\oplus O \oplus {\Cal O} \supset
\dots  \supset 
\Cal M
\oplus \Cal M\oplus {\Cal O}
$} 
\vskip -6pt 
\hbox{$\qquad\qquad\Mattt 
\qquad\quad$}  
\vskip 6pt
\hbox
{$ =
   \Cal M\oplus \Cal M\oplus {\Cal O}
\supset \dots \supset \Cal M\oplus \Cal M\oplus O \supset \dots
     \supset \Cal M\oplus \Cal M\oplus \Cal M 
\supset \dots
$} 
\vskip -6pt 
\hbox{$\qquad\qquad\Mattt 
\qquad\quad$}  
} 
}
\bigskip
\endRoster
\Roster
\Item{}  
and can also be  restored from it.
 Here the modules  isomorphic to  
${\Cal O} \oplus {\Cal O} \oplus {\Cal O}$
     do not belong to this chain  and are inserted as in  
the proof of
     Proposition~1 of \cite{P2}.
 
   \Item{(iv)} If there is an $L_{\alpha} \simeq  O \oplus O \oplus
K$ then one can restore the chain from
$$\dots \supset O \oplus O \oplus K \supset
 M \oplus O \oplus K \supset  M \oplus  M \oplus K \supset \dots
$$
     
  \Item{(v)} If there is an $L_{\alpha} \simeq  O \oplus K \oplus
K$ then the chain can be written down as 
$$\dots \supset M^i \oplus K \oplus K \supset M^{i+1} \oplus K \oplus K
\supset \dots 
$$
\endRoster
\endth

\comment 
\th Proposition 2

 Let $\{L_{\alpha} \}$ be a maximal chain of $O$-submodules in
the space $V$. There are exactly five types of maximal chains:
 \Roster   
\Item{(i)} if the chain contains a module $L$ isomorphic to $O \oplus
O \oplus O$ then all the modules of the chain are
of that type and the chain is uniquely determined by its segment 
$$\dots \supset O \oplus O \oplus O \supset
 M \oplus O \oplus O \supset
 M \oplus M \oplus O \supset
 M \oplus M \oplus M \supset
\dots $$

   \Item{(ii)} if the chain contains a module $L$ isomorphic to $O \oplus O 
\oplus {\Cal O}$ then the chain can be restored from
the segment: 
$$
\aligned 
{\Cal O} \oplus {\Cal O} \oplus {\Cal O} \supset
\dots \supset O \oplus O \oplus {\Cal O} \supset
M \oplus O \oplus {\Cal O} \supset
M \oplus M \oplus {\Cal O} \supset \dots \supset
\Cal M \oplus \Cal M\oplus {\Cal O}
} \\
 &=     
\Cal M \oplus \Cal M \oplus {\Cal O} \supset
     \dots \supset \Cal M\oplus \Cal M\oplus O \supset \Cal M\oplus
     \Cal M\oplus M \supset \dots \supset \Cal M\oplus \Cal M
\oplus
\Cal M\supset \dots 
\endaligned
$$ 
     Here the modules  isomorphic to  
${\Cal O} \oplus {\Cal O} \oplus {\Cal O}$
     do not belong to this chain  and are inserted as in the proof of
     Proposition~1 of \cite{P2}.
    
    \Item{(iii)}  all the modules $L_{\alpha} \simeq  O \oplus {\Cal O} \oplus
{\Cal  O}$. Then the chain contains a piece
$$\aligned
&
\dots \supset  {\Cal O} \oplus {\Cal O} \oplus {\Cal O} \supset
\dots \supset O \oplus {\Cal O} \oplus {\Cal O} \supset
M \oplus {\Cal O} \oplus {\Cal O} \supset \dots
 \supset
\Cal M
\oplus {\Cal O} \oplus {\Cal O} 
\\
 &=  \Cal M\oplus {\Cal O} \oplus {\Cal O}
\supset 
     \dots \supset \Cal M\oplus O \oplus {\Cal O} \supset
\dots  \supset 
\Cal M
\oplus \Cal M\oplus {\Cal O}
 \\
 &=
     \Cal M\oplus \Cal M\oplus {\Cal O}
\supset \dots \supset \Cal M\oplus \Cal M\oplus O \supset \dots
     \supset \Cal M\oplus \Cal M\oplus \Cal M  
\supset \dots
\endaligned 
     $$ 
and can also be  restored from it. 

Here the modules  isomorphic to  
${\Cal O} \oplus {\Cal O} \oplus {\Cal O}$
     do not belong to this chain  and are inserted as in the proof of
     Proposition~1 of \cite{P2}.

   \Item{(iv)} if there is an $L_{\alpha} \simeq  O \oplus O \oplus
K$ then one can restore the chain from
$$\dots \supset O \oplus O \oplus K \supset
 M \oplus O \oplus K \supset  M \oplus  M \oplus K \supset \dots
$$
     
  \Item{(v)} if there is an $L_{\alpha} \simeq  O \oplus K \oplus
K$ then the chain can be written down as 
$$\dots \supset M^i \oplus K \oplus K \supset M^{i+1} \oplus K \oplus K
\supset \dots 
$$
\endRoster
\endth 
\endcomment 

We see that the chains of the first three types correspond to two-simplices,
of  the  fourth  type  --- to  edges  of  the external boundary and the last type
represents a vertex  of the external boundary.  
As above we can glue 
     the building from apartments. To introduce them we can again repeat
     the corresponding definition for the building over a local field of
     dimension one (see Definition~4 of the previous subsection). Then the apartment $\Sigma_{.}$
     is  a union 
$$ \Sigma_{.} = \Sigma_{.}[2] \cup \Sigma_{.}[1] \cup \Sigma_{.}[0] $$
     where the pieces $ \Sigma_{.}[i]$ contain the lattices of the
     types from  $\Delta_{.}[i]$. 

The combinatorial structure of the
     apartment can be seen from two pictures at the end of the subsection. 
    There we removed the external boundary $\Sigma_{.}[0]$ which is
     simplicially isomorphic to the external boundary of an apartment of
     the building  $\Delta_{.}(PGL(3), K/K_1/k)$. The dots in the first
     picture show
     a convergence of the vertices inside the apartment. As a result the
     building is a simplicial topological space.

     We  can also describe the relations of the building with
     buildings of the same group $G$ over the complete discrete valuation  fields $K$ and
     $K_1$. In the first case there is a projection map
     $$ \pi \colon \Delta_{.}(G, K/K_1/k) \rightarrow \Delta_{.}(G, K/K_1).
     $$ 
     Under this map the big triangles containing the simplices of type (i)
     are contracted into  points, the triangles containing the simplices
     of type (ii)  go to edges and the simplices of type (iii) are
     mapped isomorphically to simplices in the target space.  The external
     boundary don't change. 
     
     The lines
\bigskip
{\ninepoint
\unitlength.9pt
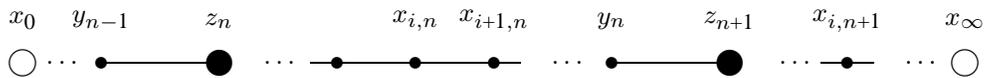
$$
\picture(360,10)

\put(0,0){\circle{10}}
\put(30,0){\circle*{5}}
\put(75,0){\circle*{10}}
\put(120,0){\circle*{5}}
\put(150,0){\circle*{5}}
\put(180,0){\circle*{5}}
\put(225,0){\circle*{5}}
\put(270,0){\circle*{10}}
\put(315,0){\circle*{5}}
\put(360,0){\circle{10}}

\put(30,0){\line(1,0){45}}
\put(110,0){\line(1,0){80}}
\put(225,0){\line(1,0){45}}
\put(305,0){\line(1,0){20}}

\put(15,0){\makebox(0,0)[cc]{$\dots$}}
\put(98,0){\makebox(0,0)[cc]{$\dots$}}
\put(208,0){\makebox(0,0)[cc]{$\dots$}}
\put(295,0){\makebox(0,0)[cc]{$\dots$}}
\put(343,0){\makebox(0,0)[cc]{$\dots$}}

\put(360,16){\makebox(0,0)[cc]{$x_{\infty}$}}
\put(75,16){\makebox(0,0)[cc]{$z_{n}$}}
\put(30,16){\makebox(0,0)[cc]{$y_{n-1}$}}
\put(150,16){\makebox(0,0)[cc]{$x_{i,n}$}}
\put(180,16){\makebox(0,0)[cc]{$x_{i+1,n}$}}
\put(270,16){\makebox(0,0)[cc]{$z_{n+1}$}}
\put(225,16){\makebox(0,0)[cc]{$y_{n}$}}
\put(315,16){\makebox(0,0)[cc]{$x_{i,n+1}$}}
\put(0,16){\makebox(0,0)[cc]{$x_{0}$}}
\endpicture
$$}\bigskip
\noindent 
     can  easily be visualized inside the apartment. Only the big white dots
     corresponding to the external boundary are missing. We have three types
     of  lines going from the inner points under the angle $2\pi/3$. They
     correspond to different embeddings of $PGL(2)$ into $PGL(3)$.
     
     Using the lines we can understand the action of the Weyl group $W$
     on an apartment. The subgroup $S_3$ acts in the same way as in
     3.2. The free subgroup  $E$ (see 3.1) has six types of
     translations along
     these three directions.  Along each line we have two opportunities
     which were introduced for $PGL(2)$.
     
Namely, if $w \in  \Gamma_K\simeq  \Bbb Z \oplus \Bbb Z \subset W$ then
$w = (0, 1)$ acts as a shift of the whole structure to the right:
$ w(x_{i,n}) = x_{i,n+2},~w(y_{n}) = y_{n+2},~w(z_{n}) = z_{n+2},
w(x_{0}) = x_{0},~w(x_{\infty}) = x_{\infty}. $
     
The element $w = (1, 0)$ acts as a shift on the points $x_{i,n}$ but leaves
fixed the points in the inner boundary
$ w(x_{i,n}) = x_{i+2,n},~w(y_{n}) = y_{n},~w(z_{n}) = z_{n},
~w(x_{0}) = x_{0},~w(x_{\infty}) = x_{\infty},  $
(see \cite{P2, The\-o\-rem~5,~v}). 

\bigskip

{\ninepoint
\unitlength.9pt
\font\edas=cmss10 scaled900
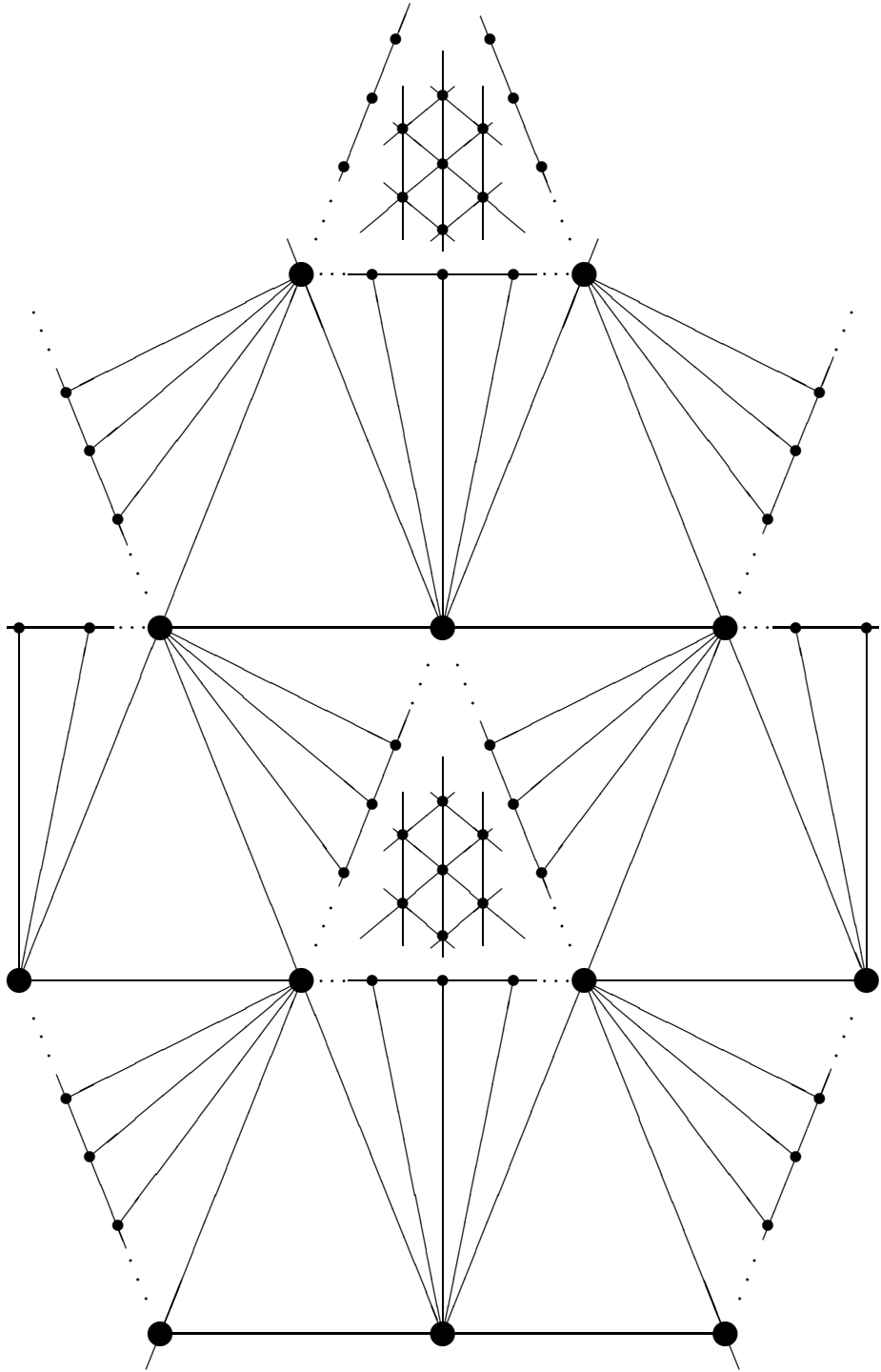
$$\picture(375,160)
%simplices type i)
\put(18,-30){{$\text{\edas simplices of type (i)}$}}
     
\put(43,5){\line(0,1){65}}
\put(60,0){\line(0,1){85}}
\put(77,5){\line(0,1){65}}

\put(25,8){\line(6,5){56}}
\put(95,8){\line(-6,5){56}}

\put(55,4){\line(6,5){30}}
\put(65,4){\line(-6,5){30}}

\put(55,70){\line(6,-5){30}}
\put(65,70){\line(-6,-5){30}}

\put(43,23){\circle*{5}}
\put(77,23){\circle*{5}}

\put(43,52){\circle*{5}}
\put(77,52){\circle*{5}}

\put(60,9){\circle*{5}}
\put(60,27){\circle*{5}}
\put(60,66){\circle*{5}}

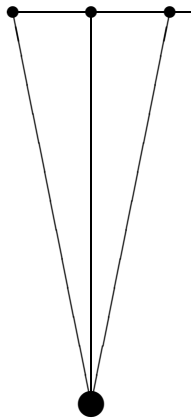
     \put(140,-30){{$\text{\edas simplices of type (ii)}$}}
     
     \put(180,0){\circle*{10}}
     \put(180,0){\line(0,1){150}}
     \put(180,0){\line(-1,5){30}}
     \put(180,0){\line(1,5){30}}
     \put(148,150){\line(1,0){70}}
     \put(150,150){\circle*{5}}
     \put(180,150){\circle*{5}}
     \put(210,150){\circle*{5}}

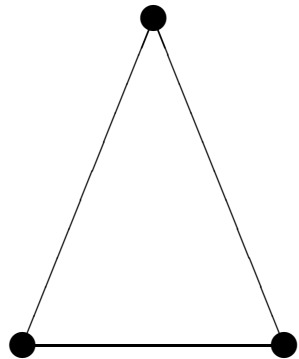
     \put(252,-30){{$\text{\edas simplices of type (iii)}$}}
     
     \put(250,0){\circle*{10}}
     \put(350,0){\circle*{10}}
     \put(300,125){\circle*{10}}
     \put(250,0){\line(1,0){100}}
     \put(250,0){\line(2,5){50}}
     \put(350,0){\line(-2,5){50}}

\endpicture
$$}

\vfill
\eject

{\ninepoint
\unitlength.85pt
$$\picture(360,540)

%\put(0,0){\line(0,1){500}}
%\put(0,0){\line(1,0){400}}

%\put(0,0){\line(1,0){400}}
%\put(0,150){\line(1,0){400}}
%\put(0,300){\line(1,0){400}}
%\put(0,450){\line(1,0){400}}
%\put(0,600){\line(1,0){400}}

%\put(0,0){\line(0,1){500}}
%\put(120,0){\line(0,1){500}}
%\put(240,0){\line(0,1){500}}
%\put(360,0){\line(0,1){500}}
%\put(480,0){\line(0,1){500}}

%triangulos
%/
\put(0,150){\line(2,5){120}}
%/
\put(60,0){\line(2,5){60}}
\put(54,-15){\line(2,5){15}}
%\
\put(60,300){\line(2,-5){120}}
%/
\put(180,0){\line(2,5){120}}
%\
\put(240,150){\line(2,-5){60}}
\put(306,-15){\line(-2,5){15}}
%\
\put(240,450){\line(2,-5){120}}
%\
\put(120,450){\line(2,-5){60}}
%/
\put(180,300){\line(2,5){60}}
\put(114,465){\line(2,-5){15}}
\put(246,465){\line(-2,-5){15}}
%
% -
\put(60,0){\line(1,0){240}}
\put(0,150){\line(1,0){120}}
\put(240,150){\line(1,0){120}}
\put(60,300){\line(1,0){240}}

%
%|
\put(0,150){\line(0,1){150}}
\put(180,0){\line(0,1){150}}
\put(180,300){\line(0,1){150}}
\put(360,150){\line(0,1){150}}
%
%nodos principales
\put(60,0){\circle*{10}}
\put(180,0){\circle*{10}}
\put(300,0){\circle*{10}}

\put(0,150){\circle*{10}}
\put(120,150){\circle*{10}}
\put(240,150){\circle*{10}}
\put(360,150){\circle*{10}}

\put(60,300){\circle*{10}}
\put(180,300){\circle*{10}}
\put(300,300){\circle*{10}}

\put(120,450){\circle*{10}}
\put(240,450){\circle*{10}}
%segmentos
%auxiliar \
%\put(0,150){\line(2,-5){60}}
%fin auxiliar \
\put(16,110){\line(2,-5){29.5}}

%\put(180,300){\line(2,-5){60}}
\put(196,260){\line(2,-5){30}}

%\put(0,450){\line(2,-5){60}}
\put(16,410){\line(2,-5){30}}
%

%auxiliar /
%\put(300,0){\line(2,5){60}}
\put(316,40){\line(2,5){29}}

%\put(120,150){\line(2,5){60}}
\put(136,190){\line(2,5){30}}

%\put(300,300){\line(2,5){60}}
\put(316,340){\line(2,5){30}}
%

%rayos |
\put(180,0){\line(1,5){30}}
\put(180,0){\line(-1,5){30}}

\put(0,150){\line(1,5){30}}

\put(360,150){\line(-1,5){30}}

\put(180,300){\line(1,5){30}}
\put(180,300){\line(-1,5){30}}
%

%segmentos -
\put(140,150){\line(1,0){80}}
\put(140,450){\line(1,0){80}}

\put(-5,300){\line(1,0){45}}
\put(320,300){\line(1,0){45}}

\put(125,149.5){$\ldots$}
\put(220,149.5){$\ldots$}
\put(40,299.5){$\ldots$}
\put(305,299.5){$\ldots$}
\put(125,449.5){$\ldots$}
\put(220,449.5){$\ldots$}
%
%rayos  ///
\put(120,150){\line(-2,-1){100}}
\put(120,150){\line(-6,-5){90}}
\put(120,150){\line(-3,-4){78}}

\put(300,300){\line(-2,-1){100}}
\put(300,300){\line(-6,-5){90}}
\put(300,300){\line(-3,-4){78}}

\put(120,450){\line(-2,-1){100}}
\put(120,450){\line(-6,-5){90}}
\put(120,450){\line(-3,-4){78}}
%
%rayos  \\\
\put(240,150){\line(2,-1){100}}
\put(240,150){\line(6,-5){90}}
\put(240,150){\line(3,-4){78}}

\put(60,300){\line(2,-1){100}}
\put(60,300){\line(6,-5){90}}
\put(60,300){\line(3,-4){78}}

\put(240,450){\line(2,-1){100}}
\put(240,450){\line(6,-5){90}}
\put(240,450){\line(3,-4){78}}

%puntos menores

\put(20,100){\circle*{5}}
\put(30,75){\circle*{5}}
\put(42,46){\circle*{5}}

\put(340,100){\circle*{5}}
\put(330,75){\circle*{5}}
\put(318,46){\circle*{5}}

\put(150,150){\circle*{5}}
\put(180,150){\circle*{5}}
\put(210,150){\circle*{5}}

\put(0,300){\circle*{5}}
\put(30,300){\circle*{5}}

\put(138,196){\circle*{5}}
\put(150,225){\circle*{5}}
\put(160,250){\circle*{5}}

\put(222,196){\circle*{5}}
\put(210,225){\circle*{5}}
\put(200,250){\circle*{5}}

\put(330,300){\circle*{5}}
\put(360,300){\circle*{5}}

\put(20,400){\circle*{5}}
\put(30,375){\circle*{5}}
\put(42,346){\circle*{5}}

\put(340,400){\circle*{5}}
\put(330,375){\circle*{5}}
\put(318,346){\circle*{5}}

\put(150,450){\circle*{5}}
\put(180,450){\circle*{5}}
\put(210,450){\circle*{5}}

%reja de arriba

\put(163,465){\line(0,1){65}}
\put(180,460){\line(0,1){85}}
\put(197,465){\line(0,1){65}}

\put(145,468){\line(6,5){56}}
\put(215,468){\line(-6,5){56}}

\put(175,464){\line(6,5){30}}
\put(185,464){\line(-6,5){30}}

\put(175,530){\line(6,-5){30}}
\put(185,530){\line(-6,-5){30}}

\put(163,483){\circle*{5}}
\put(197,483){\circle*{5}}

\put(163,512){\circle*{5}}
\put(197,512){\circle*{5}}

\put(180,469){\circle*{5}}
\put(180,497){\circle*{5}}
\put(180,526){\circle*{5}}
%
%puntitos
%\put(186.3,584){\circle*{1}}
%\put(189.5,576){\circle*{1}}
%\put(192.8,568){\circle*{1}}

\put(234,465){\circle*{1}}
\put(230.6,473){\circle*{1}}
\put(227.5,481){\circle*{1}}

\put(126.15,465){\circle*{1}}
\put(129.4,473){\circle*{1}}
\put(132.5,481){\circle*{1}}
%
%puntos medianos
\put(138,496){\circle*{5}}
\put(150,525){\circle*{5}}
\put(160,550){\circle*{5}}

\put(222,496){\circle*{5}}
\put(210,525){\circle*{5}}
\put(200,550){\circle*{5}}
\put(136,490){\line(2,5){30}}
\put(196,560){\line(2,-5){30}}
%fin
%reja de abajo

\put(163,165){\line(0,1){65}}
\put(180,160){\line(0,1){85}}
\put(197,165){\line(0,1){65}}

\put(145,168){\line(6,5){56}}
\put(215,168){\line(-6,5){56}}

\put(175,164){\line(6,5){30}}
\put(185,164){\line(-6,5){30}}

\put(175,230){\line(6,-5){30}}
\put(185,230){\line(-6,-5){30}}

\put(163,183){\circle*{5}}
\put(197,183){\circle*{5}}

\put(163,212){\circle*{5}}
\put(197,212){\circle*{5}}

\put(180,169){\circle*{5}}
\put(180,197){\circle*{5}}
\put(180,226){\circle*{5}}

%

%puntitos parte de abajo
\put(6.3,134){\circle*{1}}
\put(9.5,126){\circle*{1}}
\put(12.8,118){\circle*{1}}

\put(54,15){\circle*{1}}
\put(50.6,23){\circle*{1}}
\put(47.5,31){\circle*{1}}
\put(353.7,134){\circle*{1}}
\put(350.5,126){\circle*{1}}
\put(347.2,118){\circle*{1}}

\put(306.15,15){\circle*{1}}
\put(309.4,23){\circle*{1}}
\put(312.5,31){\circle*{1}}
%
%puntitos parte del medio 
\put(186.3,284){\circle*{1}}
\put(189.5,276){\circle*{1}}
\put(192.8,268){\circle*{1}}

\put(234,165){\circle*{1}}
\put(230.6,173){\circle*{1}}
\put(227.5,181){\circle*{1}}
\put(173.7,284){\circle*{1}}
\put(170.5,276){\circle*{1}}
\put(167.2,268){\circle*{1}}

\put(126.15,165){\circle*{1}}
\put(129.4,173){\circle*{1}}
\put(132.5,181){\circle*{1}}
%
%puntitos parte de arriba
\put(6.3,434){\circle*{1}}
\put(9.5,426){\circle*{1}}
\put(12.8,418){\circle*{1}}

\put(54,315){\circle*{1}}
\put(50.6,323){\circle*{1}}
\put(47.5,331){\circle*{1}}
\put(353.7,434){\circle*{1}}
\put(350.5,426){\circle*{1}}
\put(347.2,418){\circle*{1}}

\put(306.15,315){\circle*{1}}
\put(309.4,323){\circle*{1}}
\put(312.5,331){\circle*{1}}

\endpicture$$
}

\Bib References

\rf{BT1} 
F.  Bruhat  and J. Tits,
  {Groupes r\'eductives sur un corps local. I. 
      Donn\'ees radicielles valu\'ees}, Publ. Math. IHES {41}(1972),
      5--251; {II.Sch\'emas  en  groupes,   Existence   d'une   donn\'ee   radicielle 
        valu\'ee}, Publ. Math. IHES {60}(1984), 5--184.

\rf{BT2}
 F.  Bruhat  and J. Tits, {Sch\'emas en groupes et immeubles des groupes
                  classiques  sur  un  corps  local},  Bull.  Soc.  Math. 
                  France {112}(1984), 259--301.

\rf{FP}
T. Fimmel  and A. N.  Parshin, {An Introduction  into  the  Higher  Adelic 
Theory} (in preparation).

\rf{GK}
V. A. Ginzburg  and M. M. Kapranov, {Hecke algebras for $p$-adique  loop groups}, preprint, 1996.

\rf{M}
G. A. Mustafin, {On non-archimedean uniformization}, Math. Sbornik, 
{105}(1978), 207--237.

\rf{P1}
 A. N. Parshin, {Higher Bruhat--Tits buildings and vector bundles on an
                       algebraic surface}, Algebra and Number Theory 
        (Proc. Conf.  Inst. Exp. Math. Univ. Essen, 
                       1992), de Gruyter, Berlin, 1994, 165--192. 

\rf{P2}
A. N. Parshin, {Vector bundles and arithmetical groups I}, Proc. 
Steklov Math. Institute, {208}(1996), 212--233.

\rf{R}
M.  Ronan, {Buildings: main ideas and applications. I}, Bull.  London
          Math.Soc. {24}(1992), 1--51; {II}, ibid  
          {24}(1992), 97--126.
     
\rf{T1}
J. Tits, {On buildings and their applications}, Proc. Intern. Congr.
             Math. (Vancouver 1974), Canad. Math. Congr., Montreal, 1975,
        Vol. {1}, 209--220.

\rf{T2}
J. Tits, {Reductive groups over local fields}, Automorphic Forms,
Representations and $L$-functions, Proc. Symp. Pure Math., {33}, 
part 1, AMS, Providence, 1979, 29--70. 

\endBib

\Coordinates

Department of  Algebra,
Steklov Mathematical Institute,

Ul. Gubkina 8, Moscow GSP-1, 117966 Russia

E-mail: an\@parshin.mian.su, parshin\@mi.ras.ru 
\endCoordinates

\vfill
\pagebreak 

\end